\documentclass[11pt]{amsart}
\usepackage{extarrows}
\usepackage[colorlinks, citecolor=blue]{hyperref}
\usepackage{amsthm,amsmath,amssymb}
\usepackage{mathrsfs}
\usepackage{bold-extra}
\usepackage{algpseudocode}
\usepackage[T1]{fontenc}
 \usepackage[ruled,vlined,linesnumbered]{algorithm2e}
\usepackage{tikz}
\usepackage{tikz-cd}
\usepackage{booktabs,longtable,array}

\newtheorem{theorem}{Theorem}[section]
\newtheorem{lemma}[theorem]{Lemma}
\newtheorem{mainthm}{Theorem}

\theoremstyle{definition}
\newtheorem{definition}[theorem]{Definition}
\newtheorem{question}[theorem]{Question}

\newtheorem{proposition}[theorem]{Proposition}
\newtheorem{corollary}[theorem]{Corollary}
\newtheorem{remark}[theorem]{Remark}
\newtheorem{conjecture}[theorem]{Conjecture}

\theoremstyle{remark}

\usepackage[backref=true,backend=biber, style=alphabetic,sorting=nyt,maxnames=99,maxalphanames=99]{biblatex}
\AtBeginBibliography{\scriptsize}

\numberwithin{equation}{section}

\DeclareMathOperator{\Span}{span}
\DeclareMathOperator{\rank}{rank}

\DeclareMathOperator{\Gr}{Gr}

\DeclareMathOperator{\character}{char}
\DeclareMathOperator{\Frac}{Frac}
\DeclareMathOperator{\Spec}{Spec}
\DeclareMathOperator{\Par}{Par}

\begin{document}

\title{Cohomology Endomorphisms and the Fixed Point Property of Flag Manifolds}

\makeatletter
\def\shorttitle{Cohomology Endomorphisms and the Fixed Point Property of Flag Manifolds}
\makeatother

\author{Wangyang Lin}

\address{
\begin{tabular}{@{}l@{}}
Fudan University, Shanghai, China, \\
Queen Mary University of London, London, UK
\end{tabular}
}

\email{wylin23@m.fudan.edu.cn, wangyang.lin@qmul.ac.uk, wylin\_math@outlook.com}

%\thanks{ }

%    General info

%\date{January 1, 2001 and, in revised form, June 22, 2001.}

%\dedicatory{}

\begin{abstract}
    We study graded endomorphisms of cohomology rings and the fixed point property (FPP) of flag manifolds. Using the deformation to quantum cohomology and its curious duality, we prove Hoffman's conjecture on the cohomology rings of complex Grassmannians. More generally, we show that, uniformly over all complex Grassmannians, every graded endomorphism annihilating the first Chern class of the tautological bundle annihilates all positive-degree classes if and only if the coefficient field has characteristic zero and $-1$ is not a square.  As an application, we establish O'Neill's conjecture on the FPP of real, complex and quaternionic Grassmannians. Furthermore, we completely classify the flag manifolds of real semisimple Lie groups of exceptional type that have the FPP, including their finite products. We also show that no nontrivial flag manifold associated with a real simple Lie group whose restricted root system is of type $B$, $C$, $D$, or $BC$ has the FPP.
\end{abstract}

\maketitle

% \tableofcontents

\section{Introduction}
A topological space $X$ is said to have the \emph{fixed point property (FPP)} if every continuous self-map on $X$ has a fixed point. In O’Neill’s doctoral dissertation \cite{ONeill1974}, he investigates the fixed point property of real, complex, and quaternionic Grassmannian and raises two conjectures. Denote the Grassmannian by $\Gr(k,K^N)$, where the division ring $ \mathbb{K}=\mathbb{R},\mathbb{C},\mathbb{H}$, and let $n=N-k$. Denote the finite field of $q$ elements by $\mathbb{F}_q$.
\begin{conjecture}[{\cite[Conjecture~0.1]{ONeill1974}}]\label{conj:ONeill_FPP} 
For $\mathbb{K}=\mathbb{R},\mathbb{C}$, $\Gr(k,\mathbb{K}^N)$ has the FPP if and only if $k\ne n$ and $kn$ is even. $\Gr(k,\mathbb{H}^N)$ has the FPP if and only if $k\ne n$. 
\end{conjecture}
\begin{conjecture}[{\cite[Conjecture~0.2]{ONeill1974}}]\label{conj:ONeil_ring}
Suppose that $k\ne n$.
Every graded endomorphism $\varphi$ of $H^*(\Gr(k,\mathbb{K}^N);\mathbb{Q})$ is an Adams map, for $\mathbb{K}=\mathbb{C},\mathbb{H}$. Every graded endomorphism $\varphi$ of $H^*(\Gr(k,\mathbb{R}^N);\mathbb{F}_2)$ is an Adams map. 
\end{conjecture}
Here an Adams map acts as multiplication by $m^i$ on $H^{di}$,
where $d=2$ in the complex case and $d=4$ in the quaternionic case,
for a single $m\in\mathbb{Q}$. In the real case, it acts as multiplication
by $m^i$ on $H^i$ for a single $m\in\mathbb F_2$.
By Lefschetz's Fixed Point Theorem, Conjecture~\ref{conj:ONeil_ring} implies the if direction of Conjecture~\ref{conj:ONeill_FPP} for complex Grassmannian and even dimensional quaternionic Grassmannian.
The authors in \cite{GloverHomer1978} \cite{MilicevicRadovanovic2022} give counterexamples for real Grassmannian in Conjecture~\ref{conj:ONeil_ring}.  Hoffman classifies all the endomorphisms of $H^*(\Gr(k,\mathbb{C}^N))$ which proves Conjecture~\ref{conj:ONeil_ring} for all non-trivial endomorphisms of $H^*(\Gr(k,\mathbb{C}^N),\mathbb{Q})$ and leaves the trivial endomorphism as a conjecture, see \cite[Theorem~1.1, Conjecture~1.2]{Hoffman_Conj}.
In this paper, we prove Hoffman‘s conjecture for general coefficient fields. To state our main result, we denote by $\mathsf{Hoff}(R)$ the statement 
\begin{center}
   For all $k,n$, if a graded $R$-algebra endomorphism $\varphi$ of $H^*(\Gr(k,\mathbb{C}^N),R)$ annihilates $c_1$, then it annihilates all positive degrees.
\end{center}
\begin{mainthm}\label{mainthm:field_Hoffman}
    Let $K$ be a field. Then the following are equivalent:
    \begin{enumerate}
        \item $\mathsf{Hoff}(K)$.
        \item $\character K=0$ and $-1\notin K^2$.
    \end{enumerate}
\end{mainthm}
Since there is a ring isomorphism
\[
    H^{2*}(\Gr(k,\mathbb{C}^N);\mathbb{Q})\cong H^{4*}(\Gr(k,\mathbb{H}^N);\mathbb{Q}),
\]
see \cite[pp.\,1147]{SankaranSarkar2009}, and with Hoffman's classification \cite[Theorem~1.1]{Hoffman_Conj},
we prove Conjecture~\ref{conj:ONeil_ring} for complex and quaternionic Grassmannians. With our main theorem, we can prove Conjecture~\ref{conj:ONeill_FPP} for all real, complex and quaternionic Grassmannians.
\begin{corollary}\label{cor:FPP_CH}
    A real/complex Grassmannian has the FPP if and only if $k\ne n$ and $kn$ is even. A quaternionic Grassmannian has the FPP if and only if $k\ne n$.
\end{corollary}

The proof of Theorem~\ref{mainthm:field_Hoffman} uses a free deformation of the cohomology ring, which is called the small quantum cohomology of a Grassmannian whose special fiber is the ordinary cohomology. The special fiber is nonreduced while the quantum cohomology ring is reduced if the coefficient ring is reduced and $N$ is not a zero divisor. Then we can split each generic geometric fiber into distinct points. Specifically, on the unit fiber, we have a curious duality which is compatible with the complex conjugation when the coefficient is $\mathbb{C}$. This duality will deduce a contradiction on the unit fiber by considering hook-type Schur polynomials. 

Generally, we can raise a question on the FPP for any other flag manifold $G/P$: Which $G/P$ has the FPP? 
For real Lie groups, our flag manifolds is not the generalized flag manifolds.
\begin{mainthm}\label{mainthm:real_exc_FPP}
    Let $G$ be a connected real semisimple Lie group and let $P$ be a parabolic subgroup of $G$. Suppose that all simple Lie groups $G^{(i)}$ in the decomposition induced by the decomposition of $\mathfrak{g}$ by simple ideals
    \begin{equation}\label{eq:decomposition_flag}
        G/P \cong \prod_i \left( G^{(i)}/P^{(i)}\right),\quad P^{(i)}\subsetneq G^{(i)}
    \end{equation}
    are of real exceptional type, see \cite[Theorem~6.105]{Knapp2002}.
    Then $G/P$ has the FPP if and only if it is diffeomorphic to
   \[
(\mathbb OP^2)^a
\times
\left(\Gr(2,\mathbb{H}^4)/\langle\tau\rangle\right)^b
\times
(E_6/P_1)^c
\]
for some nonnegative integers $a,b,c$, where
\[
\tau(V)=V^\perp
\]
denotes the involution on the quaternionic Grassmannian $\Gr(2,\mathbb{H}^4)$. $P_I$ is the parabolic subgroup defined by deleting node set $I$ in the Dynkin diagram\footnote{We number the simple roots according to Bourbaki's convention.} of the restricted root system of $G$ and denote that $P_i:=P_{\{i\}}$. 
\end{mainthm}
\begin{remark}
    The three factors arise, respectively, from the real exceptional Lie groups $E\mathrm{IV}$, $E\mathrm{I}$, and $E_6$ regarded as a real Lie group. Some authors call them $E_{6(-26)}$, $E_{6(6)}$ and $E_6(\mathbb{C})$.
\end{remark}
As a corollary, we can immediately identify the FPP for complex semisimple Lie groups, since the only complex Lie group components in Theorem~\ref{mainthm:real_exc_FPP} is $E_6/P_1$.
\begin{corollary}\label{cor:exc_FPP}
    Let $G$ be a connected complex semisimple Lie group whose simple factors of its Lie algebra are all of exceptional type, and let $P$ be a proper parabolic subgroup of $G$. Then $G/P$ has the FPP if and only if $G/P$ is isomorphic to a product of $E_6/P_1$. 
\end{corollary}
\begin{mainthm}\label{mainthm:BCD_FPP}
    Suppose that $G$ is a real semisimple Lie group and $P$ is a parabolic subgroup. If one of the simple Lie group in decompositino~\eqref{eq:decomposition_flag} has restricted root system of $B_n$, $C_m$, $D_r$ or $BC_p$ type, then $G/P$ does not have the FPP.
\end{mainthm}
\begin{corollary}
    Let $G$ be a connected complex semisimple Lie group whose simple factors of its Lie algebra are all of $B_n$, $C_m$ or $D_{r}$ type, and let $P$ be a proper parabolic subgroup of $G$. Then $G/P$ does not have the FPP. 
\end{corollary}

For the FPP of general mainfolds, the story becomes complicated.  A product of closed manifold $X=\prod_i X_i$ may not have the FPP even if each $X_i$ has the FPP. A counterexample with different $X_i$'s is constructed in \cite{Husseini1977}. If we addtionally suppose that the factors are the same, the FPP problem is still open.
\begin{question}[\cite{KwasikSun2017}]
    Suppose that a closed manifold $X$ has the FPP, does $X\times X$ have the FPP?
\end{question}
\begin{question}\label{que:FPP_flag_manifolds}
    Suppose that $X,Y$ are in one of these categories: homogeneous spaces, homogeneous spaces of compact Lie groups, flag manifolds of semisimple real/complex Lie groups and they have the FPP. Does $X\times X$ or $X\times Y$ have the FPP? 
\end{question}
See Appendix~\ref{sec:append_flag_FPP} for a partial answer for Question~\ref{que:FPP_flag_manifolds} for flag manifolds.

\subsection*{Organization of this paper} From Section~\ref{sec:Pre} to \ref{sec:mainthm}, we prove our main theorem for complex Grassmannians. In Section~\ref{sec:Pre}, we review some basic properties of a Grassmannian and a deformation family of its cohomology ring, which is exactly the small quantum cohomology of a Grassmannian. In Section~\ref{sec:curious}, we introduce a curious duality on the unit fiber of the deformation family. In Section~\ref{sec:mainthm}, we finish the proof of our main theorem. In Section~\ref{sec:FPP}, we focus on the FPP and prove Corollary~\ref{cor:FPP_CH} on real/complex/quaternionic Grassmannian and Theorem~\ref{mainthm:real_exc_FPP}, \ref{mainthm:BCD_FPP} flag manifolds associated to $B,C,D$ and exceptional type Lie groups.

\subsection*{Acknowledgment} The author would like to thank Alex Fink, Xinrui Wang, Andy Xiaoan Yang for helpful discussions concerning this problem, and Ping Li for bringing this problem to his attention several years ago.

\subsection*{Disclosure of AI use} The author used GPT-6-Astra during the early stages of this work to explore possible proof strategies. GPT-6-Astra identified the key proof idea of using curious duality on the unit fiber to restrict the Schur expansion to hook partitions. The author developed the geometric interpretation of the argument, supplied the details needed to complete the proof, and determined exactly which coefficient fields satisfy the vanishing statement for all complex Grassmannians. 
All mathematical arguments were reviewed and verified by the author, who takes full responsibility for the content of this paper.

\section{Preliminaries}\label{sec:Pre}
\subsection{Quantum Schubert calculus} 
Let $R$ be a commutative ring with identity. Throughout, ring homomorphisms preserve the identity. Let
\[
    \Lambda_R=R[e_1,\ldots,e_k]=R[x_1,\ldots,x_k]^{S_k}
\]
be the symmetric polynomial ring in $k$ variables, where $e_i$ is the $i$-th elementary symmetric polynomial. We use polynomial degree, so that $\deg e_i=i$, rather than cohomological degree. For a graded ring or module $A$, we write $A_d$ for its degree-$d$ component and $A_{>0}$ for its positive-degree part. 
To identify the relation of the cohomology ring, we need the following generating function and its inverse
\[
    E(z)=1+\sum_{i=1}^k e_i z^i,\quad E(z)^{-1}= 1+\sum_{j=1}^\infty (-1)^j h_j z^j,
\]
where $h_j$ is the $j$-th complete symmetric polynomial. Set $e_0=h_0=1$. By Schubert calculus, see for example \cite{Fulton_young_tab}, we have the following isomorphism for a complex Grassmannian
\[
A_R:=H^{2*}(\Gr(k,\mathbb{C}^N);R)\cong \Lambda_R/I_R,\quad I_R=(h_{n+1},\ldots,h_N),\quad N=k+n.
\]
In $A_R$, we have $e_i=0$ for $i>k$.
We can assume that $k\le n$ by Grassmannian duality.
Write
\[
\begin{aligned}
    \Par(k,n)=\{\alpha=(\alpha_1\ge\cdots\ge\alpha_k\ge0):\alpha_1\le n\}.
\end{aligned}
\]
We denote Schur polynomial by $s_\alpha$. $\{s_\alpha:\alpha\in\Par(k,n)\}$ forms an $R$-basis of $A_R$. With some abuse of notation, we denote both those symmetric polynomials in $\Lambda$ and their images under quotient maps/base changes by $e_i, h_j, s_\alpha$.

Denote the last relation $h_N$ by $(-1)^{k-1}Q$, we define
\[
B_R:=R[Q,e_1,\ldots,e_k]/
(h_{n+1},\ldots,h_{N-1},h_N-(-1)^{k-1}Q),\quad \deg(Q)=N,
\]
which is a deformation family over $R[Q]$ with special fiber $A_R$ at $0$. In fact, $B_R$ is exactly the \emph{small quantum cohomology} of $\Gr(k,\mathbb{C}^N)$, see \cite[\S\,8]{Buch_quantum}. We denote that
\[
    J_R=(h_{n+1},\ldots,h_{N-1}),
\]
then $B_R=\Lambda_R/J_R$ and $A_R=B_R/(Q)$.

\begin{proposition}\label{prop:free}
For every commutative ring $R$,
\begin{equation}\label{eq:free}
B_R=\bigoplus_{\alpha\in\Par(k,n)}R[Q]s_\alpha,
\qquad
(B_R)_m=\bigoplus_{\substack{\alpha\in\Par(k,n),\ a\ge0\\|\alpha|+aN=m}}R Q^a s_\alpha.
\end{equation}
For $b\in R$, denote the fiber by $B_{R,b}=B_R/(Q-b)$; its Schur polynomials form an $R$-basis.
\end{proposition}

\begin{proof}
    In \cite[\S\,2]{Buch_quantum}, the author proves the $\mathbb{Z}[Q]$-module isomorphism $B_\mathbb{Z}\cong A_\mathbb{Z}\otimes_\mathbb{Z}\mathbb{Z}[Q]$.
    Then we obtain first the decomposition by tensoring $R$. And the second decomposition comes from decomposition of $R[Q]$ by degree in $Q$. It is preserved by base change since $-\otimes_RR'$ distributes over direct sums. 

    For each fiber $B_{R,b}$, we can view $R$ as an $R[Q]$-algebra by $Q\to b$, then the decomposition is obtained by using base change $B_{R,b}=B_R\otimes_{R[Q]}R$. 
\end{proof}
\begin{remark}
    The geometric meaning of this proposition can be interpreted as a free deformation family $\Spec B_R\to \mathbb{A}_R^1$ with special fiber $\Spec A_R$ at zero. In fact, $\Spec A_R$ is a fat point with length $\binom{N}{k}$, which is nonreduced. The nilpotent elements obstruct the proof inside $A_R$, so we need to split the fat point into $\binom{N}{k}$ distinct points which is one of the generic fiber.  
\end{remark}

\subsection{Deformation family of cohomology} By Proposition~\ref{prop:free}, we know that each fiber $B_{R,b}$ is spanned by $s_\alpha$ with $\alpha\in\Par(k,n)$ over coefficient ring $R$. And the geometric point of $\Spec B_{R,b}$ is $\binom{N}{k}$ distinct points. We first identify each geometric point of $\Spec B_{R,b}$ as a pair of generating functions.

Before specializing on each fiber, we prove the following presentation of $B_R$.
\begin{lemma}[{\cite[pp.\,479 (6)]{Postnikov_quantum}}]\label{lemma:repre_2}
There is another presentation of the small quantum cohomology of $\Gr(k,\mathbb{C}^N)$, that is 
    \begin{equation}\label{eq:repre_2}
        B_R \cong R[e_1,\ldots,e_k,h_1,\ldots,h_n,Q]/([z^m](E(z)T(z)-1-(-1)^nQz^N):1\le m\le N),
    \end{equation}
    where
    \[
        T(z)=\sum_{j=0}^n(-1)^jh_jz^j.
    \]
\end{lemma}

The relations in this presentation is called \emph{quantum Witten relations}, see \cite{Witten1995Verlinde}.
With the special representation in Lemma~\ref{lemma:repre_2}, we can describe each geometric fiber of $B_R$ over $R[Q]$. Suppose that $L$ is an algebraically closed field with $\character L\nmid N$ and let $s:\Spec L\to \Spec R[Q]$ be a geometric point satisfying $b:=s^*(Q)\in L^\times$. Then the geometric fiber at $s$ is determined by the base change of $\Spec B_R\to \mathbb{A}_R^1$ to $\Spec L$, that is the spectrum of $B_R\otimes_{R[Q],s^*}L$. Moreover, each generic geometric fiber is reduced.
The description of this geometric fiber and a more general reducedness of $B_R$ are in the following proposition.
\begin{proposition}\label{prop:reduced}
    \begin{enumerate}
        \item Suppose that $L$ is an algebraically closed field with $\character L\nmid N$ and let $s:\Spec L\to \Spec R[Q]$ be a geometric point satisfying $b:=s^*(Q)\in L^\times$. There is an $L$-algebra isomorphism 
        \[
            B_R\otimes_{R[Q],s^*}L\cong \prod_{S\subset S_b,|S|=k} L,\quad S_b=\{x:x^N=(-1)^{k-1}b\},
        \]
        with the multiplication in the product of $L$ defined as multiplications in components.
        % \item $\Spec B_R[(NQ)^{-1}]$ is finite and \'{e}tale over $\Spec R[N^{-1},Q,Q^{-1}]$ of rank $\binom{N}{k}$. 
        \item If $R$ is reduced and $N$ is not a zero divisor in $R$, then $B_R$ is reduced. In particular, if $K$ is a field with $\character K\nmid N$, then $B_K$ is reduced.
    \end{enumerate}
\end{proposition}

\begin{proof}
    (1). Set $d=\binom{N}{k}$. By Proposition~\ref{prop:free}, $B_R$ is a free $R[Q]$-module of rank $d$. The geometric fiber $C=B_R\otimes_{R[Q],s^*}L$ is an $L$-algebra of dimension $d$.

    For some $L$-point $p$ of $C$, that is an $L$-algebra homomorphism $p^*:C\to L$, let
    \[
        E_p(z)=1+\sum_{i=1}^k p^*(e_i)z^i,\quad T_p(z) =  \sum_{j=0}^n (-1)^j p^*(h_j)z^j,
    \]
    and define
    \[
        P(u)=u^kE_p(-u^{-1}),\quad U(u)=u^nT_p(-u^{-1}).
    \]
    The relations in Lemma~\ref{lemma:repre_2} are equivalent to
    \begin{equation}\label{eq:P-U}
        P(u)U(u)=u^N-(-1)^{k-1}b.
    \end{equation}
    Since $\character L\nmid N$ and $b\in L^\times$, this equation has $N$ distinct roots.

    With a subset $S$ containing $k$ roots, define $P=\prod_{x\in S}(u-x)$ and $U=\prod_{S_b\backslash S}(u-x)$. Then we can recover $E_p$ and $T_p$ which define the geometric point $p^*:C\to L$. Conversely, every geometric point $p^*$ gives the decomposition in \eqref{eq:P-U}, and $P$ has $k$ roots. Thus, geometric points $C\to L$ are in bijection with the $k$-subset of $N$ roots. 

    For each geometric point $p_i^*$ of $C$, $\ker p^*_i$ is a maximal ideal since it is an $L$-algebra homomorphism. Then by Chinese Remainder Theorem, the map $(p_i^*)_{i=1}^d:C\to L^d$ is surjective. This is an isomorphism by comparing the dimensions.

    % (2). $B_R[(NQ)^{-1}]$ is a finite free $R[N^{-1},Q,Q^{-1}]$-algebra with finitely many generators and relations. Thus, it is flat and finitely presented. By (1), we know the geometric fiber is $L^d$ and it is \'{e}tale over $L$. Finally, it is an \'etale morphism.

    \noindent (2). Suppose that $K$ is a field and $\character K\nmid N$. Freeness of $B_K$ guarantees
    \[
        B_K\hookrightarrow B_K\otimes_{K[Q]}K(Q)\hookrightarrow B_K\otimes_{K[Q]}\overline{K(Q)}\cong\overline{K(Q)}^d.
    \]
    The last isomorphism comes from (1). Because $B_K$ embeds into a reduced ring, it is reduced.

    For general $R$ which is reduced and $N$ is not a zero divisor, we have
    \[
        R\hookrightarrow R[N^{-1}]\hookrightarrow\prod_{\mathfrak{p}\in\Spec R[N^{-1}]}\Frac(R[N^{-1}]/\mathfrak{p}),
    \]
    where $R[N^{-1}]$ is reduced. The last inclusion holds since the kernel of the map is the nilradical.

    In $B_R[N^{-1}]$, we can write each element as an $R[N^{-1},Q]$-linear combination of $s_\alpha$.
    \[
        f=\sum_{a,\alpha} c_{a,\alpha} Q^as_\alpha,\quad c_{a,\alpha}\in R[N^{-1}].
    \]
    We define the inclusion
    \[
        B_{R[N^{-1}]}\hookrightarrow \prod_{\mathfrak{p}\in\Spec R[N^{-1}]} B_{\Frac(R[N^{-1}]/\mathfrak{p})},
    \]
    by sending each $c_{a,\alpha}$ into each residue field $\Frac(R[N^{-1}]/\mathfrak{p})$. 
    This map is injective since the kernel of the diagonal map is
    \[
        \ker\left(R[N^{-1}]\to\prod_{\mathfrak{p}\in\Spec R[N^{-1}]} \Frac(R[N^{-1}]/\mathfrak{p})\right)=\bigcap_\mathfrak{p}\mathfrak{p}=\sqrt{(0)}=(0).
    \]
    Since the characteristic of each residue field does not divide $N$, the target is reduced. Thus, $B_R$ is reduced by $B_R\hookrightarrow B_{R[N^{-1}]}$.
\end{proof}

In Proposition~\ref{prop:reduced}, we prove the reducedness for each geometric fiber and for $B_R$ when $R$ is reduced when $N$ is not a zero divisor. We will see in the next proposition how reducedness plays a role in our proof.
\begin{proposition}\label{prop:extension}
    Suppose that $k\le n$. For every graded $R$-algebra endomorphism $\varphi$ on $A_R$, there is a unique $F_i\in (\Lambda_R)_i$ whose image under $\Lambda_R\to A_R$ is $\varphi(e_i)$. This defines  $\Phi:\Lambda_R\to \Lambda_R$ and induces a graded $R$-algebra endomorphism of $B_R$ which we denote by $\Phi_B$ such that
    \[
        \Phi_B(Q)=\lambda Q, \quad \lambda\in R.
    \]
    Moreover, if $B_R$ is reduced, then $\lambda=0$ is equivalent to vanishing of $\varphi$ on $(A_R)_{>0}$.
\end{proposition}
\begin{proof}
    Denote the quotient map by $\pi:\Lambda_R\to A_R$. Since $\pi$ gives $(\Lambda_R)_i\cong (A_R)_i$ on each degree $0\le i\le n$, there is a unique $F_i\in(\Lambda_R)_i$ such that $\pi(F_i)=\varphi(e_i)$. This gives a unique map $\Phi:\Lambda_R\to\Lambda_R$ such that $e_i\mapsto F_i$ which preserves degree. Since $\pi\Phi=\varphi\pi$, for any $P\in I_R$, we have $\pi\Phi(P)=\varphi\pi(P)=0$ which means $\Phi(I_R)\subset I_R$.
    Then, for $h_j$ with $n<j<N$, we have
    \[
        \Phi(h_j) = \sum_{t=n+1}^j a_{j-t} h_t,\quad a_{j-t}\in (\Lambda_R)_{j-t},
    \]
    which means $h_N$ is not in the summation and then $\Phi(J_R)\subset J_R$. Thus, $\Phi_B$ is a ring homomorphism well-defined on $B_R$. And by $(I_R)_N = (J_R)_N +Rh_N$, we have $\Phi(h_N)-\lambda h_N\in J_R$ for some $\lambda\in R$, that is $\Phi_B(Q)=\lambda Q$. And the $\lambda$ is uniquely determined since $B_R$ is a free $R[Q]$-module.

    Finally, suppose that $B_R$ is reduced. If $\lambda=0$, then $\Phi_B$ factors through $A_R$ and  $(A_R)_{>0}$ is sent to the nilradical of $B_R$ which is zero. Thus, $\varphi((A_R)_{>0})=0$ by modulo $Q$. Conversely, if $\varphi$ vanishes on $(A_R)_{>0}$, then $F_i=0$ for every $1\le i\le k$. $Q=(-1)^{k-1}h_N$ can be represented as a polynomial of $e_1,\ldots,e_k,h_{n+1},\ldots,h_{N-1}$, then $\Phi_B(Q)=0$.
\end{proof}

\section{Curious Duality}\label{sec:curious}

\subsection{Curious duality}
In \cite{Postnikov_quantum}, an involution (called \emph{curious duality}) on the localization $B_R[Q^{-1}]$ is defined to relate Gromov–Witten invariants of complementary degrees and to determine the highest quantum degree in a Schubert product from the lowest one. This duality is used in our proof to prove a hook type spanning proposition later. Let us first introduce the definition of this duality.
\begin{definition}
    For $\alpha\in\Par(k,n)$, define the \emph{Durfee number} to be $r(\alpha)=\#\{i:\alpha_i\ge i\}$. When $r(\alpha)=1$, the partition is of \emph{hook-type}. Define an involution $\star$ on $\Par(k,n)$ preserving the Durfee number as
    \[
        \alpha^\star_i=\begin{cases}
            n+r(\alpha)-\alpha_{r(\alpha)+1-i},\quad i\le r(\alpha),\\
            r(\alpha)-\alpha_{k+r(\alpha)+1-i},\quad i>r(\alpha),
        \end{cases}
    \]
    and $|\alpha^\star|=Nr(\alpha)-|\alpha|$.
\end{definition}

The diagram operation keeps the Durfee square and complements and rotates each of the right and lower arms inside its own rectangle. For example, for $k=4,n=6$ it sends $(6,4,1,0)$ to $(4,2,2,1)$; both have Durfee number $2$, and their sizes add to $20=2(4+6)$.
\begin{center}
\definecolor{starsquare}{RGB}{105,142,195}
\definecolor{starright}{RGB}{105,189,210}
\definecolor{starlower}{RGB}{237,177,105}
\begin{tikzpicture}[x=0.65cm,y=0.65cm]
\node at (3,0.65) {$\alpha=(6,4,1,0)$};
\filldraw[fill=starsquare,draw=gray!65,line width=0.3pt] (0,0) rectangle (1,-1);
\filldraw[fill=starsquare,draw=gray!65,line width=0.3pt] (1,0) rectangle (2,-1);
\filldraw[fill=starright,draw=gray!65,line width=0.3pt] (2,0) rectangle (3,-1);
\filldraw[fill=starright,draw=gray!65,line width=0.3pt] (3,0) rectangle (4,-1);
\filldraw[fill=starright,draw=gray!65,line width=0.3pt] (4,0) rectangle (5,-1);
\filldraw[fill=starright,draw=gray!65,line width=0.3pt] (5,0) rectangle (6,-1);
\filldraw[fill=starsquare,draw=gray!65,line width=0.3pt] (0,-1) rectangle (1,-2);
\filldraw[fill=starsquare,draw=gray!65,line width=0.3pt] (1,-1) rectangle (2,-2);
\filldraw[fill=starright,draw=gray!65,line width=0.3pt] (2,-1) rectangle (3,-2);
\filldraw[fill=starright,draw=gray!65,line width=0.3pt] (3,-1) rectangle (4,-2);
\filldraw[fill=white,draw=gray!65,line width=0.3pt] (4,-1) rectangle (5,-2);
\filldraw[fill=white,draw=gray!65,line width=0.3pt] (5,-1) rectangle (6,-2);
\filldraw[fill=starlower,draw=gray!65,line width=0.3pt] (0,-2) rectangle (1,-3);
\filldraw[fill=white,draw=gray!65,line width=0.3pt] (1,-2) rectangle (2,-3);
\filldraw[fill=gray!12,draw=gray!65,line width=0.3pt] (2,-2) rectangle (3,-3);
\filldraw[fill=gray!12,draw=gray!65,line width=0.3pt] (3,-2) rectangle (4,-3);
\filldraw[fill=gray!12,draw=gray!65,line width=0.3pt] (4,-2) rectangle (5,-3);
\filldraw[fill=gray!12,draw=gray!65,line width=0.3pt] (5,-2) rectangle (6,-3);
\filldraw[fill=white,draw=gray!65,line width=0.3pt] (0,-3) rectangle (1,-4);
\filldraw[fill=white,draw=gray!65,line width=0.3pt] (1,-3) rectangle (2,-4);
\filldraw[fill=gray!12,draw=gray!65,line width=0.3pt] (2,-3) rectangle (3,-4);
\filldraw[fill=gray!12,draw=gray!65,line width=0.3pt] (3,-3) rectangle (4,-4);
\filldraw[fill=gray!12,draw=gray!65,line width=0.3pt] (4,-3) rectangle (5,-4);
\filldraw[fill=gray!12,draw=gray!65,line width=0.3pt] (5,-3) rectangle (6,-4);
\draw[line width=1pt] (0,0) rectangle (2,-2);
\node at (12,0.65) {$\alpha^\star=(4,2,2,1)$};
\filldraw[fill=starsquare,draw=gray!65,line width=0.3pt] (9,0) rectangle (10,-1);
\filldraw[fill=starsquare,draw=gray!65,line width=0.3pt] (10,0) rectangle (11,-1);
\filldraw[fill=starright,draw=gray!65,line width=0.3pt] (11,0) rectangle (12,-1);
\filldraw[fill=starright,draw=gray!65,line width=0.3pt] (12,0) rectangle (13,-1);
\filldraw[fill=white,draw=gray!65,line width=0.3pt] (13,0) rectangle (14,-1);
\filldraw[fill=white,draw=gray!65,line width=0.3pt] (14,0) rectangle (15,-1);
\filldraw[fill=starsquare,draw=gray!65,line width=0.3pt] (9,-1) rectangle (10,-2);
\filldraw[fill=starsquare,draw=gray!65,line width=0.3pt] (10,-1) rectangle (11,-2);
\filldraw[fill=white,draw=gray!65,line width=0.3pt] (11,-1) rectangle (12,-2);
\filldraw[fill=white,draw=gray!65,line width=0.3pt] (12,-1) rectangle (13,-2);
\filldraw[fill=white,draw=gray!65,line width=0.3pt] (13,-1) rectangle (14,-2);
\filldraw[fill=white,draw=gray!65,line width=0.3pt] (14,-1) rectangle (15,-2);
\filldraw[fill=starlower,draw=gray!65,line width=0.3pt] (9,-2) rectangle (10,-3);
\filldraw[fill=starlower,draw=gray!65,line width=0.3pt] (10,-2) rectangle (11,-3);
\filldraw[fill=gray!12,draw=gray!65,line width=0.3pt] (11,-2) rectangle (12,-3);
\filldraw[fill=gray!12,draw=gray!65,line width=0.3pt] (12,-2) rectangle (13,-3);
\filldraw[fill=gray!12,draw=gray!65,line width=0.3pt] (13,-2) rectangle (14,-3);
\filldraw[fill=gray!12,draw=gray!65,line width=0.3pt] (14,-2) rectangle (15,-3);
\filldraw[fill=starlower,draw=gray!65,line width=0.3pt] (9,-3) rectangle (10,-4);
\filldraw[fill=white,draw=gray!65,line width=0.3pt] (10,-3) rectangle (11,-4);
\filldraw[fill=gray!12,draw=gray!65,line width=0.3pt] (11,-3) rectangle (12,-4);
\filldraw[fill=gray!12,draw=gray!65,line width=0.3pt] (12,-3) rectangle (13,-4);
\filldraw[fill=gray!12,draw=gray!65,line width=0.3pt] (13,-3) rectangle (14,-4);
\filldraw[fill=gray!12,draw=gray!65,line width=0.3pt] (14,-3) rectangle (15,-4);
\draw[line width=1pt] (9,0) rectangle (11,-2);
\draw[->,thick] (6.8,-2) -- (8.2,-2);
\node at (7.5,-1.35) {$\star$};
\end{tikzpicture}
\end{center}

With the $\star$ involution on partitions, we have an involution on the localized quantum cohomology ring.
\begin{theorem}[\emph{curious duality}, {\cite[Theorem~6.5]{Postnikov_quantum}}]
    The involution $\tilde{D}:Q\mapsto Q^{-1}$ defines an involution  on $B_R[Q^{-1}]$ by
    \[
        \tilde{D}(s_\alpha)=Q^{-r(\alpha)}s_{\alpha^\star}.
    \]
\end{theorem}

\begin{corollary}\label{cor:inversion}
The curious duality induces an involution on $B_{R,1}$ (which we also denote by $\tilde{D}$), where the factor $Q^{-r(\alpha)}$ disappears. Its coordinate formulas are
\begin{equation}\label{eq:invcoords}
\tilde{D}(e_i)=e_{k-i}h_n,\quad
\tilde{D}(h_j)=e_kh_{n-j}.
\end{equation}
\end{corollary}

\begin{proof}
    Since the curious duality maps $Q-1$ to $-Q^{-1}(Q-1)$, which means it keeps the ideal $(Q-1)$, then it induces an involution on $B_{R,1}$ such that $\tilde{D}(s_\alpha)=s_{\alpha^\star}$.

    For $1\le i\le k$,
    \[
        \tilde{D}(e_i)=\tilde{D}(s_{(1^i)})=s_{(1^i)^\star}=s_{(n,1^{k-i})}.
    \]
    When $k=i$, $\tilde{D}(e_i)=s_{(n)}=e_0h_n$. When $k>i$,
    \[
        \tilde{D}(e_i)=s_{(n,1^{k-i})}=h_ne_{k-i}-s_{(n+1,1^{k-i-1})}=h_ne_{k-i}-\sum_{t=0}^{k-i-1}(-1)^th_{n+1+t}e_{k-i-1-t}=h_ne_{k-i}.
    \]
    The last equality holds since the subscripts of $h_j$ in the summation range from $n+1$ to $N-i$.

    And for $1\le j \le n$, similarly,
    \[
        \tilde{D}(h_j)=\tilde{D}(s_{(j)})=s_{(j)^\star}=s_{(n+1-j,1^{k-1})}=e_kh_{n-j}. \qedhere
    \]
    \let\qed\relax
\end{proof}

\subsection{Curious duality and conjugation over $\mathbb{C}$}\label{subsec:curious_for_c}

\begin{lemma}[{\cite[Theorem~2.12]{hengelbrock2002involutionquantumcohomologyring}}]\label{lemma:evaluation}
    Let $p\in\Spec B_{\mathbb{C},1}$ be a complex point defined by a $\mathbb{C}$-algebra homomorphism
    \[
        p^*:B_{\mathbb{C},1}\to \mathbb{C},
    \]
    and let $\{x_1,\ldots,x_k\}$ be its corresponding set of roots in Proposition~\ref{prop:reduced}.
    Then any Schur polynomial $s_\alpha$ in $B_{\mathbb{C},1}$ which can be viewed as a regular function on $\Spec B_{\mathbb{C},1}$ has its evaluation at each point as
    \[
        s_\alpha(p)=p^*(s_\alpha)=s_\alpha(x_1,\ldots,x_k),\quad x_i\in \mathbb{C}.
    \]
    Moreover,
    \[
        \tilde{D}s_\alpha(p) = \overline{s_\alpha(p)}.
    \]
\end{lemma}
\begin{proof}
    For each point $p\in \Spec B_{\mathbb{C},1}$, the corresponding polynomials $E_p$ and $T_p$ as shown in the proof of Proposition~\ref{prop:reduced} satisfy
\[
    E_p(z)T_p(z)=1+(-1)^nz^N.
\]
And this corresponds to $P(u)$ and $U(u)$ such that
\[
    P(u)U(u)=\prod_{x\in S}(u-x)\prod_{x\notin S}(u-x)=u^N-(-1)^{k-1}.
\]
Thus each point $p$ is in bijection with
\[
    S=\{x_1,\ldots,x_k\}\subset\{x:x^N=(-1)^{k-1}\}\text{ in }\mathbb{C}.
\]
For each $x$, the inverse $x^{-1}=\overline{x}$. Each $e_i$ is a regular function on $\Spec B_{\mathbb{C},1}$, then
\[
    e_i(p)=p^*(e_i)=e_i(x_1,\ldots,x_k)
\]
by the definition of $E_p$ and $P$.

By Corollary~\ref{cor:inversion} and $e_kh_n=1$ we have
\[
    e_i(\tilde{D}^*(p))=\tilde{D}e_i(p)=\frac{e_{k-i}(x_1,\ldots,x_k)}{e_k(x_1,\ldots,x_k)}=e_i(x_1^{-1},\ldots,x_k^{-1})=\overline{e_i(p)}.
\]
By Giambelli's formula, we finish the proof.
\end{proof}

We have proved that the curious duality has a nice property on $B_{\mathbb{C},1}$. Then we need to know about its interaction with $\Phi$ defined in Proposition~\ref{prop:extension}. However, $\Phi_B$ can not induce a endomorphism on the unit fiber, so we need to scale it as following.

Suppose that $\lambda$ in Proposition~\ref{prop:extension} is not zero for $\Phi_B$ on $B_\mathbb{C}$. Set $a\in \mathbb{C}$ such that $a^N=\lambda$ and define a $\mathbb{C}$-algebra endomorphism on $\Lambda_\mathbb{C}$ by setting
\[
    \Psi(e_i) = a^{-i}\Phi(e_i)=a^{-i}F_i.
\]
This also induces a $\Psi_B:B_\mathbb{C}\to B_\mathbb{C}$ such that $\Psi_B(Q)=Q$, since $Q=(-1)^{k-1}h_N$ is a homogeneous symmetric polynomial of degree $N$. Moreover, it induces a ring endomorphism $\psi:B_{\mathbb{C},1}\to B_{\mathbb{C},1}$. We define 
\[
    \sigma: \sum c_{\alpha} s_\alpha \mapsto \sum \overline{c_\alpha}s_\alpha, 
\]
then we define $\overline{\Psi}=\sigma\Psi\sigma$ which induces $\overline{\psi}:B_{\mathbb{C},1}\to B_{\mathbb{C},1}$.

\begin{lemma}\label{lemma:twist}
    On the unit fiber $B_{\mathbb{C},1}$, we have
    \[
        \tilde{D}\psi=\overline{\psi}\tilde{D},
    \]
    where we restrict $\tilde{D}$ on $B_{\mathbb{C},1}$.
\end{lemma}

\begin{proof}
Suppose that
\[
    g = \sum_\alpha c_\alpha s_\alpha,\quad \psi(s_\alpha)=\sum_{\beta}a_{\alpha\beta}s_\beta,\quad c_\alpha, a_{\alpha\beta}\in\mathbb{C}.
\]
The maps $\tilde{D},\psi,\overline{\psi}$ are all $\mathbb{C}$-linear. For each complex point $p\in\Spec B_{\mathbb{C},1}$, the left-hand side is
\[
    \tilde{D}\psi(g)(p) = \tilde{D}\left(\sum_{\alpha,\beta} c_\alpha a_{\alpha\beta}s_\beta\right)(p) = \sum_{\alpha,\beta} c_\alpha a_{\alpha\beta} \tilde{D}(s_\beta)(p)=\sum_{\alpha,\beta} c_\alpha a_{\alpha\beta}\overline{s_\beta(p)}.
\]
whereas the right-hand side is
\[
     \overline{\psi}\tilde{D}(g)(p) = \sigma \psi \sigma \left(\sum_\alpha c_\alpha \tilde{D}(s_\alpha)\right)(p) =  \sum_{\alpha}c_\alpha\sigma\psi(s_{\alpha^\star})(p).
\]
We prove a necessary equality
\[
    \sigma \psi (s_{\alpha^\star})(p)=\overline{\sigma\psi(s_\alpha)(p)},
\]
by applying the argument of Lemma~\ref{lemma:evaluation} for $b_i=\sigma\psi(e_i)(p)$. The polynomials
\[
    E_b(z)=1+\sum_{i=1}^k b_iz^i,\quad T_b(z)=\sum_{j=0}^n (-1)^jh_j(b_1,\ldots,b_k)z^j
\]
satisfy $E_bT_b=1+(-1)^nz^N$. Then $b_i=e_i(y_1,\ldots,y_k)$ with $y_i\in \mathbb{C}$ being the roots of equation $u^N=(-1)^{k-1}$. Since
\[
    \tilde{D}(e_i)=e_{k-i}h_n,\quad e_kh_n=1,
\]
then
\[
    \sigma\psi\tilde{D}(e_i)(p)=\frac{e_{k-i}(y_1,\ldots,y_k)}{e_k(y_1,\ldots,y_k)}=e_i(y_1^{-1},\ldots,y_k^{-1})=\overline{\sigma\psi(e_i)(p)}.
\]
And by Giambelli's formula we prove the equality. Then we finish the whole proof by
\[
    \overline{\psi}\tilde{D}(g)(p)=\sum_{\alpha}c_\alpha\sigma\psi(s_{\alpha^\star})(p)=\sum_{\alpha}c_\alpha\overline{\sigma\psi(s_\alpha)(p)}=\sum_{\alpha,\beta}c_\alpha\overline{\sigma a_{\alpha\beta}s_\beta(p)}=\sum_{\alpha,\beta} c_\alpha a_{\alpha\beta}\overline{s_\beta(p)}.
\]
Since $B_{\mathbb C,1}$ is reduced, equality at all complex points implies equality in $B_{\mathbb C,1}$.
\end{proof}

\section{Proof of Theorem~\ref{mainthm:field_Hoffman}}\label{sec:mainthm}

In this section, we prove Theorem~\ref{mainthm:field_Hoffman} by showing
\[
    (\neg \mathsf{Hoff}(K)) \wedge (\character K=0) \Rightarrow -1\in K^2,
\]
and
\[
    (\character K>0) \vee (-1\in K^2) \Rightarrow \neg\mathsf{Hoff}(K).
\]
For the first implication, $\neg \mathsf{Hoff}(K)$ implies that there is some $\Gr(k,\mathbb{C}^N)$ on whose cohomology $\varphi:A_K\to A_K$ is non-trivial on positive degrees. By Proposition~\ref{prop:extension}, the extension $\Phi_B$ on $B_K$ satisfies $\Phi_B(Q)=\lambda Q$ with $\lambda\ne 0$. Together with the curious duality, we will obtain some requirements on $K$.

For the second implication, we display two examples in $\Gr(2,\mathbb{C}^4)$ and $\Gr(2,\mathbb{C}^{p^2-1}),p\ge 3$ which violate $\mathsf{Hoff}(K)$.

\subsection{Hook-type support} Recall that a partition $\alpha$ or a Schur polynomial $s_\alpha$ is of \emph{hook-type} if its Durfee number $r(\alpha)=1$. We call a symmetric polynomial \emph{$K$-hook-spanned} if it is a $K$-linear combination of hook-type Schur polynomials.

\begin{proposition}\label{prop:hook_span}
    Suppose that $K$ is a field of characteristic zero and $2\le k\le n$. If the extension $\Phi:\Lambda_K\to \Lambda_K$ in Proposition~\ref{prop:extension} satisfies $\lambda\ne 0$ for $\Phi_B$, then $\Phi(e_i)$ and $\Phi(h_j)$ are both $K$-hook-spanned for $1\le i\le k$ and $1\le j\le n$.
\end{proposition}

\begin{proof}
 \textbf{Step 1.} $K=\mathbb{C}$.

 Denote $G_i=\Psi(e_i)$ where $\Psi$ is defined in Section~\ref{subsec:curious_for_c}. Let $\pi_1$ be the projection from $\Lambda_\mathbb{C}$ to $B_{\mathbb{C},1}$, then $\pi_1(G_i)=\psi\pi_1(e_i)=\psi(e_i)$. By Corollary~\ref{cor:inversion} and Lemma~\ref{lemma:twist} we have
 \begin{equation}\label{eq:hook_equality}
    \tilde{D}\pi_1(G_i) = \pi_1 (\overline{\Psi}(e_{k-i})\overline{\Psi}(h_n)).
 \end{equation}
 We compare the sizes of the partitions in the Schur expansions. The right-hand side has weights $|\alpha|=N-i$ for all partitions in its Schur expansion since it comes from a degree $N-i$ polynomial in $B_\mathbb{C}$. And by
 \[
    \pi_1(G_i) = \sum_{\substack{\alpha\in\Par(k,n)\\|\alpha|=i}}c_\alpha s_\alpha,\quad \tilde{D}\pi_1(G_i)=\sum_{\substack{\alpha\in\Par(k,n)\\|\alpha|=i}}c_\alpha s_{\alpha^\star}.
 \]
 and $|\alpha^\star|=Nr(\alpha)-i$ where $r(\alpha)$ is the Durfee number, then we must have
 \[
    Nr(\alpha)-i = N - i,\quad r(\alpha)=1,
 \]
 for all $c_\alpha\ne 0$ to balance the weights of partitions in \eqref{eq:hook_equality}. This means $G_i$'s are $\mathbb{C}$-hook-spanned. For $H_j=\Psi(h_j)$, we have the similar equality
 \[
    \tilde{D}\pi_1(H_j) = \pi_1 (\overline{\Psi}(h_{n-j})\overline{\Psi}(e_k)).
 \]
 The same argument shows that $H_j$'s are $\mathbb{C}$-hook-spanned. By the definition of $\Psi$, we have both $\Phi(e_i)$ and $\Phi(h_j)$ are $\mathbb{C}$-hook-spanned.

 \textbf{Step 2.} Field $K$ with zero characteristic.

 Suppose that
 \[
    F_i := \Phi(e_i) = \sum_{\alpha} c_{i\alpha}s_\alpha,\quad \Phi:\Lambda_K\to\Lambda_K, \quad c_{i\alpha}\in K.
 \]
 Consider $K_0=\mathbb{Q}(c_{i,\alpha}:1\le i\le k,\alpha\in\Par(k,n),\lambda)$. This field admits an embedding into $\mathbb{C}$. So coefficients in front of non-hook-type Schur polynomial in $\Phi(e_i)$ are zero in $K_0$ then are zero in $K$ since $K_0\subset K$ for every $\Phi$. The proof is the same for $\Phi(h_j)$.
\end{proof}

\begin{lemma}\label{lemma:square_hook}
    Suppose that $K$ is a field and $2\le m\le k\le n$. If both
    \[
        0\ne F=\sum_{j=0}^{m-1}a_js_{(m-j,1^j)}\in A_K,
    \]
    and $F^2$ are $K$-hook-spanned in $A_K$, then $a_0\ne 0$ and $a_0^2+a_1^2=0$. In particular, $-1\in K^2$.
\end{lemma}

\begin{proof}
    First,
    \[
        [s_{(m,m)}]F^2=a_0^2+a_1^2=0,
    \]
    by the Littlewood–Richardson rule in $A_K$.
    
    The inclusion $\Gr(r,\mathbb{C}^{n+r})\hookrightarrow\Gr(k,\mathbb{C}^{n+k})$ induces a pullback $p_r$ of symmetric polynomials such that
    \[
        p_r (s_\alpha) =\begin{cases}
            s_\alpha,\quad \ell(\alpha)\le r,\\
            0,\quad \text{otherwise},
        \end{cases}
    \]
    and $\alpha_1$ is still no greater than $n$. In particular, the pullback of a hook-type Schur polynomial is still of hook-type. Assume that $a_0=a_1=0$, then there is a minimal $t\ge 2$ such that $a_t\ne 0$. Let $r=t+1, a=m-r$, then
    \[
        p_r(F) = a_t s_{(a+1,1^{r-1})}=a_te_rh_a,
    \]
    and by expanding $h_a^2$ and then multiplying by $e_r^2$ we have
    \[
        (p_r(F))^2= a_t^2\sum_{s=0}^a s_{(2a-s+2,s+2,2^{r-2})}.
    \]
    The partition $s_{(a+2,a+2,2^{r-2})}$ has Durfee number $2$ and its coefficient in $(p_r(F))^2$ is $a_t^2\ne0$. This partition is in $\Par(k,n)$ since $r\le k$, $a+2\le n$. This coefficient still appears in $F^2$ that makes a contradiction. So $a_0\ne 0$ and $a_1/a_0$ is a square root of $-1$.
\end{proof}

\subsection*{Proof of Theorem~\ref{mainthm:field_Hoffman}}
We finish the proof by proving two statements for $2\le k\le n$. When $k=1$, the proof is trivial on the cohomology of a complex projective space.

\textbf{Statement 1. }
\[
    (\neg \mathsf{Hoff}(K)) \wedge (\character K=0) \Rightarrow -1\in K^2,
\]

\begin{proof}
By Proposition~\ref{prop:extension}, there is some $\Gr(k,\mathbb{C}^N),k\le n$ such that  $\lambda\ne 0$ for the extension $\Phi_B:B_K\to B_K$ induced by $\Phi:\Lambda_K\to\Lambda_K$. 
Denote the quotient map by $\pi:\Lambda_K\to A_K$ with $\varphi \pi=\pi \Phi$.

Since $\varphi(e_1)=0$, there is a minimal $m\ge 2$ such that $\varphi(e_m)\ne 0$ and $\varphi(e_i)=0$ for $i>k$. Expanding the inverse of $1+\sum_{i\ge m}\varphi(e_i)z^i$ we obtain
\[
    \varphi(h_{2m}) = \varphi(e_m)^2-\varphi(e_{2m})\text{ in }A_K.
\]
When $2m\le n$, 
\[
    \varphi(e_{m})=\pi\Phi(e_{m}),\quad\varphi(e_{2m})=\pi\Phi(e_{2m}),\quad\varphi(h_{2m})=\pi\Phi (h_{2m})
\]
 are all $K$-hook-spanned by Proposition~\ref{prop:hook_span}, which means $(\pi\Phi(e_m))^2=\varphi(e_m)^2$ is $K$-hook-spanned. 
When $2m>n$, both $e_{2m}$ and $h_{2m}$ are zero and $K$-hook-spanned, then $\varphi(e_m)^2=0$ is obviously $K$-hook-spanned.
Thus, $-1\in K^2$ by Lemma~\ref{lemma:square_hook}.
\end{proof}

\textbf{Statement 2.}
\[
    (\character K>0) \vee (-1\in K^2) \Rightarrow \neg\mathsf{Hoff}(K).
\]

\begin{proof}
    \textbf{Case I.} If $\character K = 0$ and $\mathbf{i}^2=-1\in K$, then consider $A_K$ with $k=n=2$. Pieri's formula gives
\[
    h_2^2=e_2^2=s_{(2,2)},\quad h_2e_2=0.
\]
Let $F=h_2+\mathbf{i}e_2$ which is non-trivial in $A_K$ and then $F^2=0$. Set $\varphi(e_1)=0,\varphi(e_2)=F$, which is non-trivial on $(A_K)_{>0}$. It is an endomorphism since $\varphi(h_3)=\varphi(h_4)=0$ by
\[
    1+\sum_{j\ge1} \varphi(h_j)z^j=(1+Fz^2)^{-1}=1-Fz^2+F^2z^4+\cdots.
\]

\textbf{Case II.} If $\character K = 2$, $-1$ is always a square in $K$, then the example in \textbf{Case I.} is still valid.

\textbf{Case III.} If $\character K = p>2$, then consider $A_K$ with $k=2$ and $n=p^2-3$. Set $F=e_1^2-4e_2=(x-y)^2$ where $x,y$ are two Chern roots. Then the Frobenius map gives
\[
    F^{(p^2-1)/2}=(x-y)^{p^2-1}=\frac{x^{p^2}-y^{p^2}}{x-y}=h_{p^2-1}.
\]
Define
\[
    \varphi(e_1)=0,\quad \varphi(e_2)=F,
\]
which is an endomorphism since
\[
    1+\sum_{j\ge1} \varphi(h_j)z^j=(1+Fz^2)^{-1},
\]
and then $\varphi(h_{p^2-2})=0$ and $\varphi(h_{p^2-1})=(-F)^{(p^2-1)/2}=0$ in $A_K$. Since $F=s_{(2)}-3s_{(1,1)}$, $\varphi$ is non-trivial on $(A_K)_{>0}$.
\end{proof}

\subsection{Hoffman's Conjecture with arbitrary coefficients} For any reduced ring $R\ne 0$, we have the following inclusion
\[
    R\hookrightarrow \prod_{\mathfrak{p}\in\operatorname{Min}(R)} R_\mathfrak{p},
\]
where $\operatorname{Min}(R)$ is the set of minimal prime ideals of $R$. This map is injective because the localization $R_\mathfrak{p}$ is the residue field $R_\mathfrak{p}$ for every minimal prime ideal.

Before moving forward to the proof of the main theorem, let us first prove the following lemma of localization.
\begin{lemma}\label{lemma:localization}
    Suppose that $S\subset R$ is a multiplicatively closed subset. Then
    \[
        \mathsf{Hoff}(R)\Rightarrow \mathsf{Hoff}(S^{-1}R).
    \]
\end{lemma}
\begin{proof}
    Suppose that $\varphi':A_{S^{-1}R}\to A_{S^{-1}R}$ satisfies $\varphi'(e_1)=0$ but
    \[
        f_i:=\varphi'(e_i)\ne 0 \text{ in }A_{S^{-1}R}
    \]
    for some $i\ge 2$.
    Choose $u\in S$ such that there is $F_i\in (A_R)_i$ for each $i\ge 2$ such that
    \[
        F_i/1 = u^if_i\text{ in } A_{S^{-1}R},
    \]
    and set $F_1=0$. For $n<j\le N$, let $H_j=h_j(F_1,\ldots,F_k)$ in $A_R$, then
    \[
        H_j/1 = h_j (uf_1,\ldots,u^kf_k)=u^jh_j(f_1,\ldots,f_k)=0,
    \]
    by relations in $I_{S^{-1}R}$. Thus, there is some $v_j\in S$ such that $v_j H_j=0$. Let $v=\prod_{j=n+1}^Nv_j$ and $G_i=v^iF_i$, then
    \[
        h_j(G_1,\ldots,G_k)=v^jh_j(F_1,\ldots,F_k)=v^jH_j=0.
    \]
    Thus, the map
    \[
        \varphi:A_R\to A_R, e_i\mapsto G_i
    \]
    is well-defined. $\varphi(e_1)=0$ since $G_1=vF_1=0$. For a non-trivial $f_i$, in the localization
    \[
        G_i/1 = (uv)^if_i\ne 0,
    \]
    since $uv\in S$ is a unit. Then $G_i\ne 0$. $\varphi$ is non-trivial on $(A_R)_{>0}$.
\end{proof}

\begin{theorem}\label{thm:ring_Hoffman}
    Let $R\ne0$ be a commutative ring. Then the following are equivalent:
    \begin{enumerate}
        \item $\mathsf{Hoff}(R)$.
        \item $R$ is reduced and  $\mathsf{Hoff}(R_\mathfrak{p})$ for every minimal prime ideal $\mathfrak{p}$, where $R_\mathfrak{p}$ is the localization at $\mathfrak{p}$. $R_\mathfrak{p}$ is exactly the residue field $\Frac(R/\mathfrak{p})$ of $R$ at $\mathfrak{p}$ when $R$ is reduced.
        \item $R$ is torsion-free over $\mathbb{Z}$ and $a^2+b^2=0$ implies $a=b=0$ in $R$.
    \end{enumerate}
\end{theorem}
\begin{proof}
(1) $\Leftrightarrow$ (2). Assume that $x\ne 0$ is a nilpotent element in $R$ with $x^t=0,x^{t-1}\ne 0$. Then for $\Gr(2,\mathbb{C}^4)$ and
\[
    e_1\mapsto 0,\quad e_2\mapsto x^{t-1}e_2,
\]
the inverse of $1+x^{t-1}e_2z^2$ is $1-x^{t-1}e_2z^2$, which means $\varphi(h_3)=\varphi(h_4)=0$. So this is a ring endomorphism of $A_R$ non-trivial on $e_2$. Thus, $R$ is reduced. 
And by Lemma~\ref{lemma:localization}, $\neg\mathsf{Hoff}(R_\mathfrak{p})$ implies $\neg\mathsf{Hoff}(R)$.

Conversely, suppose that $R$ is reduced and $\mathsf{Hoff}(R_\mathfrak{p})$ holds for every minimal prime ideal. Let $\varphi:A_R\to A_R$ be a graded $R$-algebra endomorphism satisfying $\varphi(e_1)=0$ and we suppose that
\[
    \varphi(e_i) = \sum_{|\alpha|=i}c_{i\alpha} s_\alpha,\quad c_{i\alpha}\in R,
\]
then vanishing of $\varphi$ on $e_1$ is preserved by base change to $R_\mathfrak{p}$. Extend $\varphi$ to the endomorphism on $A_{R_\mathfrak{p}}$, 
\[
    \varphi(e_i)=\sum_{|\alpha|=i}\eta_\mathfrak{p}(c_{i\alpha}) s_\alpha,
\]
where $\eta_\mathfrak{p}:R\to R_\mathfrak{p}$ is the canonical homomorphism. By $\mathsf{Hoff}(R_\mathfrak{p})$, $\eta_\mathfrak{p}(c_{i\alpha})=0$ in $R_\mathfrak{p}$. That is 
\[
    c_{i\alpha}\in\bigcap_{\mathfrak{p}\in\operatorname{Min}(R)}\ker \eta_\mathfrak{p}=\bigcap_{\mathfrak{p}\in\operatorname{Min}(R)}\mathfrak{p}=(0),
\]
by the reducedness of $R$. Thus, $\varphi(e_i)=0$ for $1\le i\le k$.

(2) $\Leftrightarrow$ (3). 
By Theorem~\ref{mainthm:field_Hoffman}, $\mathsf{Hoff}(R_\mathfrak{p})$ is equivalent to $\character R_\mathfrak{p}=0$ and $-1\notin R_\mathfrak{p}^2$. If this holds for every minimal prime ideal, then a torsion element in $R$ must be zero in each $R_\mathfrak{p}$. Thus it must be zero in $R$ by the reducedness. If $a^2+b^2=0$, then for each $R_\mathfrak{p}$ both $a$ and $b$ should be zero, thus they are both zero in $R$ by reducedness.

Conversely, assume (3). If there is a nonzero nilpotent element $x\in R$ with minimal $t$ such that $x^t=0,x^{t-1}\ne 0$, then $a=x^{t-1},b=0$ is the nonzero solution of $a^2+b^2=0$. Thus, $R$ is reduced.

Suppose that there is some $\mathfrak{p}$ such that $R_\mathfrak{p}$ has positive characteristic $p$, then $p/1=0$ and there is some $s\notin \mathfrak{p}$ such that $ps=0$ in $R$ which violates the torsion-free condition. And if $(a/s)^2=-1$ in $R_\mathfrak{p}$ where $s\notin\mathfrak{p}$, then there is some $t\notin \mathfrak{p}$ such that $t(a^2+s^2)=0$. Then
\[
    (ta)^2+(ts)^2=t^2(a^2+s^2)=0.
\]
By the condition in (3), $ts=0$ but $t,s\notin\mathfrak{p}$ which violates the definition of a prime ideal. So $-1\notin R_\mathfrak{p}^2$. Finally, by Theorem~\ref{mainthm:field_Hoffman} (3) implies (2).
\end{proof}

\section{Fixed Point Property}\label{sec:FPP}

\subsection{Fixed Point Property for Grassmannians} 
Recall that a topological space is said to have FPP if every continuous self-map on it has a fixed point. Before our proof of our Corollary~\ref{cor:FPP_CH}, we recall the cohomology ring of a quaternionic Grassmannian. Suppose that $S$ is the quaternionic tautological bundle over $\Gr(k,\mathbb{H}^N)$ and $S_\mathbb{C}$ is its underlying complex vector bundle. Let $x_j,-x_j,j=1,\ldots,k$ be the Chern roots of $S_\mathbb{C}$
and the symplectic Pontryagin class is defined by $q_i=(-1)^ic_{2i}(S_\mathbb{C})$. Then $q_i$ is the $i$-th elementary symmetric polynomial in $x_j^2s$. In \cite[pp.\,1147]{SankaranSarkar2009}, the authors give an isomorphism 
\[
    H^{2*}(\Gr(k,\mathbb{C}^N);\mathbb{Q})\cong H^{4*}(\Gr(k,\mathbb{H}^N);\mathbb{Q}),
\]
which preserves the partition of Schubert cells. So
\[
    H^{4*}(\Gr(k,\mathbb{H}^N);\mathbb{Q})=\mathbb{Q}[q_1,\ldots,q_k]/(h_{n+1},\ldots,h_N),
\]
where $h_r$ is the $r$-th complete symmetric polynomials in $x_j^2$'s.

\subsection*{Proof of Corollary~\ref{cor:FPP_CH} for complex Grassmannians} Suppose that $k\le n$.
When $k<n$ and $kn$ is even, we can prove that every endomorphism of $A_\mathbb{Q}$ is an Adams map by Theorem~\ref{mainthm:field_Hoffman} and the classification of endomorphisms in \cite[Theorem~1.1]{Hoffman_Conj}. 
For every continuous self-map $f$, there is an $m\in\mathbb{Z}$ such that $f^*(e_i)=m^ie_i$. Then the Lefschetz number
\begin{equation}\label{eq:lef_num}
        L(f)=\sum_{i=0}^{kn} (-1)^{2i}\operatorname{tr}(f^*|_{H^{2i}})=\sum_{i=0}^{kn}b_{2i}m^i=\begin{bmatrix}N\\k\end{bmatrix}_m=\prod_{i=1}^k\frac{1-m^{n+i}}{1-m^i},
\end{equation}
where $\left[\begin{smallmatrix} N \\ k \end{smallmatrix}\right]_m$ is the Gaussian binomial polynomial. The second equality is because $f^*|_{H^{2i}}$ is scaling by $m^i$. The third equality is the Poincar\'{e} polynomial of a complex Grassmannian.

Since the roots of a Gaussian binomial polynomial are on the unit circle and $b_{2i}>0$, we have $L(f)\ne 0$ for $m\ne -1$. If $m=-1$, by $q$-Lucas theorem, when $kn$ is even
\[
    \begin{bmatrix}N\\k\end{bmatrix}_{-1} = \binom{\lfloor N/2\rfloor}{\lfloor k/2\rfloor}\ne 0.
\]
By the Lefschetz's Fixed Point Theorem, $f$ has a fixed point.

When $k=n$, let $f(V)=V^\perp$ under some Hermitian metric, and $f$ has no fixed point. When $k<n$ and $kn$ is odd, we can view $\mathbb{C}^N$ as a right $\mathbb{H}$-module, and take $J(v)=v\mathbf{j}$ such that $J^2=-\operatorname{Id}$.
Let $f(V)=J(V)$, then $J(V)=V$ implies $J^2|_V=-\operatorname{Id}$ which means $\dim_\mathbb{C}V=k$ is even, but it must be odd. So $f$ has no fixed point.\hfill$\Box$

\subsection*{Proof of Corollary~\ref{cor:FPP_CH} for quaternionic Grassmannians} Suppose that $k\le n$. When $k<n$, since there is a ring isomorphism
\[
    H^{2*}(\Gr(k,\mathbb{C}^N);\mathbb{Q})\cong H^{4*}(\Gr(k,\mathbb{H}^N);\mathbb{Q}),
\]
the Lefschetz number of every $f$ is still that in \eqref{eq:lef_num} with Betti numbers replaced by $b_{4i}$. So when $m\ne -1$, $f$ has a fixed point. And when $n=k$, $f(V)=V^\perp$ has no fixed point.

When $m=-1$ and $kn$ is even, the Gaussian binomial polynomial is still nonzero. But it is zero when $m=-1$ and $kn$ is odd. So we end the proof by showing that there is no endomorphism induced by self-map on $H^*(\Gr(k,\mathbb{H}^N))$ with $m=-1$. A proof by mod $3$ Steenrod cyclic reduced power operation for quaternionic projective space can be found in \cite[Section~IV.15]{Bredon1993}. Before extending the method to quaternionic Grassmannian, we recall the mod $p$ Steenrod cyclic reduced power operation.

\begin{definition}
For an odd prime $p$, the Steenrod cyclic reduced power operation $\mathcal{P}^k,k\ge 0$  is a natural transformation
\[
    \mathcal{P}^k:H^i(X;\mathbb{F}_p) \to H^{i+2k(p-1)}(X;\mathbb{F}_p)
\]
which is a homomorphism satisfying the following axioms:
\begin{enumerate}
    \item $\mathcal{P}^0=\operatorname{id}$.
    \item If $\deg(a)=2k$, then $\mathcal{P}^k(a)=a^p$.
    \item If $\deg(a)<2k$, then $\mathcal{P}^k(a)=0$.
    \item (Cartan formula) $\mathcal{P}^k(a\smile b)=\sum_{j=0}^k\mathcal{P}^j(a)\smile\mathcal{P}^{k-j}(b).$
\end{enumerate}
\end{definition}
For the Chern roots $x_j,-x_j,j=1,\ldots,k$ of $S_\mathbb{C}$ over $\Gr(k,\mathbb{H}^N)$,
\[
    \mathcal{P}^0(x_j)=x_j,\quad\mathcal{P}^1(x_j)=x_j^p,\quad \mathcal{P}^k(x_j)=0,k \ge 2.
\]
So for the first symplectic Pontryagin class,
\[
    \mathcal{P}^1(q_1) = \sum_{j=1}^k \mathcal{P}^1(x_j^2)=2\sum_{j=1}^k x_j^{p+1}.
\]
When $p=3$, we have $\mathcal{P}^1(q_1)=-q_1^2-q_2\mod 3$.
If $f^*(q_1)=-q_1$, then $f^*(q_2)=(-1)^2q_2=q_2$ and $f^*(\mathcal{P}^1(q_1))=-q_1^2-q_2\mod 3$ and $\mathcal{P}^1(f^*(q_1))=q_1^2+q_2\mod 3$. By naturality, we have $q_1^2+q_2=0\mod 3$.
When $k\ge 2$, $q_1^2$ and $q_2$ form a basis of $H^8(\Gr(k,\mathbb{H}^N);\mathbb{F}_3)$. When $k=1$, $q_1^2$ forms a basis of $H^8(\Gr(k,\mathbb{H}^N);\mathbb{F}_3)$. So $q_1^2+q_2\mod 3$ can not be zero. So $m\ne -1$.\hfill$\Box$

\subsection*{Proof of Corollary~\ref{cor:FPP_CH} for real Grassmannians}
For a real Grassmannian, if $k=n$, $f(V)=V^\perp$ still has no fixed point. If $k<n$ and $kn$ is odd, then a $J:\mathbb{R}^N\to\mathbb{R}^N$ with $J^2=-\operatorname{Id}$ can be still defined since $N=n+k$ is even. $f(V)=J(V)$ has no fixed point.

Conversely, suppose that $k<n$ and $kn$ is even.
Let $S$ be the tautological bundle and $Q$ is the quotient bundle. Denote the Pontryagin classes $p=p(S)$ and $\bar{p}=p(Q)$. Let $p(t)$ and $\bar{p}(t)$ be the generated function of $p_i$ and $\bar{p}_i$. When $kn$ is even and $2$ is invertible in $R$, the cohomology ring of a real Grassmannian is
\[
    H^{*}(\Gr(k,\mathbb{R}^n);R)\cong \frac{R[p_1,\ldots,p_a,\bar{p}_1,\ldots,\bar{p}_b]}{([t^d](p(t)\bar{p}(t))-1):d=1,\ldots,a+b)}\quad a=\left\lfloor\frac{k}{2}\right\rfloor,b=\left\lfloor\frac{n}{2}\right\rfloor.
\]
And there is an isomorphism
\[
    H^{4*}(\Gr(k,\mathbb{R}^N);R)\cong H^{2*}(\Gr(a,\mathbb{C}^{a+b});R),\quad p_i\mapsto e_i,\bar{p}_j\mapsto h_j.
\]
When $R=\mathbb{Q}$, suppose that every self-map $f$ satisfies $f^*(p_1)=m p_1$ for some rational number $m$.
By Theorem~\ref{mainthm:field_Hoffman} and Theorem~1.1 in \cite{Hoffman_Conj},  $f^*$ can only have two cases: (1) $f^*$ as an Adams map when $a< b$ or $a=b, m=0$, or (2) $f^*(p_i)=m^i\bar{p}_i$ for some $m\ne 0$.

For case (1), we use mod $3$ Steenrod cyclic reduced power operation, and then the argument is the same as the quaternionic Grassmannians which shows that $L(f)\ne0$. In fact, this argument does not need $a<b$.

For case (2), $a=b$ and $k<n$ imply $k=2a$ and $n=2a+1$. So we consider the Lefschetz number on
\[
    H^{4*}(\Gr(2a,\mathbb{R}^{4a+1});\mathbb{Q})\cong H^{2*}(\Gr(a,\mathbb{C}^{2a});\mathbb{Q}).
\]
On a the Schur class of $\Gr(a,\mathbb{C}^{2a+1})$, $f^*(s_\lambda)=(-m)^{|\lambda|} s_{\lambda'}$ where $\lambda'$ is the transpose. Then
\[
    L(f) = \sum_{\substack{\lambda\in\Par(a,a)\\\lambda=\lambda'}} (-m)^{|\lambda|}=\prod_{j=1}^a(1-m^{2j-1}).
\]
Thus, $L(f)\ne 0$ if and only if $m\ne 1$. When $m=1$, consider the oriented real Grassmannian $\widetilde{\Gr}(2a,\mathbb{R}^{4a+1})$ which is a $2$ to $1$ universal covering of the real Grassmannian. Denote the covering map by $\pi$ and the Euler class of the oriented Grassmannian by $e$, then $e^2=\pi^*(p_a)$. We still write $p_i,\bar{p}_i$ for their pullbacks. The cohomology of the oriented ring Grassmannian can be repesented as
\[
    A':=H^{4*}(\widetilde{\Gr}(2a,\mathbb{R}^{4a+1});\mathbb{Q})\cong \frac{\mathbb{Q}[p_1,\ldots,p_a,\bar{p}_1,\ldots,\bar{p}_a,e]}{([t^k](p(t)\bar{p}(t))-1):k=1,\ldots,2a, e^2=p_a)}.
\]
By induction,
\[
    \bar{p}_a=-p_a-\sum_{i=1}^{a-1}p_i\bar{p}_{a-1}=-p_a+G(p_1,\ldots,p_{a-1})=-e^2+G.
\]
Suppose that the lifting $\tilde{f}$ satifies $\tilde{f}^*(e)=ce+F(p_1,\ldots,p_a)$, then
\[
    (ce+F)^2=\tilde{f}^*(e)^2=\tilde{f}^*(p_a)=\bar{p}_a.
\]
This gives a relation in cohomological degree $4a$ in $A'$. 
We simplify $A'$ as
\[
    \frac{\mathbb{Q}[p_1,\ldots,p_a,e]}{(h_{a+1},\ldots,h_{2a}, e^2=p_a)}, \quad \sum_{j\ge 0}h_jt^j=p(t)^{-1}.
\]
All degrees in relations of $A'$ is at least $4a+4$. So $m\ne 1$.
By the Lefschetz Fixed Point Theorem, we finish the proof.\hfill$\Box$

\begin{remark}
    The only possible octonionic projective spaces are $\mathbb{O}P^1\cong \mathbb{S}^8,\mathbb{O}P^2$. Only $\mathbb{O}P^2$ has the FPP. And there is no general family of classical Grassmannians over the octonions.
\end{remark}

\subsection{Flag manifolds of exceptional Lie groups} For a semisimple Lie group $G$ and one of its parabolic subgroup $P$, we have the decomposition
\[
    G/P \cong \prod_i \left( G^{(i)}/P^{(i)}\right).
\]
Throughout, a parabolic subgroup is understood to be the normalizer of a parabolic subalgebra.

When $G$ is simple, and choose a Iwasawa decomposition $G=KAN$, where $K$ is the maximal compact group, $A$ is the abelian subgroup and $N$ is the nilpotent part. And the decomposition on Lie algebra is $\mathfrak{g}=\mathfrak{k}\oplus\mathfrak{a}\oplus\mathfrak{n}$ which is a direct sum of real vector space.

Suppose $\Sigma$ is the restricted root system of $G$ on $\mathfrak{a}$ and each parabolic subgroup $P=P_I$ can be defined by deleting subset $I$ from its simple roots $\Delta$. Denote that
\[
    \mathfrak{a}_I= \{H\in\mathfrak{a}:\alpha(H)=0,\alpha\in\Delta\backslash I\},\quad Z_K(\mathfrak{a}_I)=\{h\in K:\operatorname{Ad}(h)(H)=H,H\in\mathfrak{a}_I\}.
\]
Let $A_I$ be the torus of $\mathfrak{a}_I$.
By Langlands decomposition $P_I=M_IA_IN_I$ \cite[Proposition~7.82(a)~7.83(d)]{Knapp2002}, there is a diffeomorphism $G/P_I\cong K/(K\cap P_I)\cong K/Z_K(\mathfrak{a}_I)$. 
$Z_G(\mathfrak{a}_I)$ is called the standard Levi subgroup, and $Z_K(\mathfrak{a}_I)=Z_G(\mathfrak{a}_I)\cap K$\footnote{The Levi groups can be different but they are conjugate in $P$. $Z_G(\mathfrak{a}_I)$ is the standard one which is compatible with the Cartan involution. }.
With the normalizer
\[
    N_K(\mathfrak{a})=\{h\in K:\operatorname{Ad}(h)(\mathfrak{a})=\mathfrak{a}\},
\]
we can define the restricted Weyl group $W(G,\mathfrak{a})=N_K(\mathfrak{a})/Z_K(\mathfrak{a})$, which is isomorphic to the Weyl group of restricted root system \cite[Theorem~6.57]{Knapp2002}. For a complex Lie group, its restricted root system is the root system of its complex Lie algebra.
\begin{lemma}\label{lemma:weyl_FPP}
    If $-\operatorname{Id}_\mathfrak{a}$ is an element of the restricted Weyl group, any nontrivial $G/P$ does not have FPP.
\end{lemma}
\begin{proof}
    Take an representation element $\sigma$ of $-\operatorname{Id}_\mathfrak{a}$, it preserves $\mathfrak{a}_\Theta$ and then $Z_K(\mathfrak{a}_\Theta)$. It defines a self-map on $K/Z_K(\mathfrak{a}_\Theta)$ by sending $hZ_K(\mathfrak{a}_\Theta)$ to $h\sigma Z_K(\mathfrak{a}_\Theta)$, which has no fixed point. If it has a fixed point, then $\sigma\in Z_K(\mathfrak{a}_\Theta)$ but $\sigma$ action on $\mathfrak{a}_\Theta$ is not identity.
\end{proof}
With this lemma we can rule out a large amount of exceptional real simple Lie group. By \cite[pp. 82]{Humphreys1990}, only Weyl group of reduced root system $A_n(n\ge 2), D_{2m+1}, E_6$ does not contain $-\operatorname{Id}$. 
\begin{corollary}\label{cor:BCD_FPP}
        Suppose that $G$ is a real semisimple Lie group and $P$ is a parabolic subgroup. If one of the simple Lie group in decompositino~\eqref{eq:decomposition_flag} has restricted root system of $B_n$, $C_r$, $D_{2m}$ or $BC_p$ type, then $G/P$ does not have the FPP.
\end{corollary}
\begin{proof}
    The FPP on $G/P$ implies the FPP on each components. 
\end{proof}
This corollary excludes most types for Theorem~\ref{mainthm:BCD_FPP}. We will study $D_{2m+1}$-type and prove Theorem~\ref{mainthm:BCD_FPP} at the end of this section.

\begin{corollary}\label{cor:exc_remaining}
    For real simple Lie algebras whose complexifications are exceptional,
only the following three exceptional types satisfy that $-\operatorname{Id}_{\mathfrak{a}}$ is not in their Weyl groups.
\[
    E\mathrm{I},\quad E\mathrm{IV},\quad E_6\text{ as a real Lie algebra}.
\]
\end{corollary}
\begin{proof}
    Category (a) of \cite[Theorem~6.105]{Knapp2002} consists of complex Lie groups as real Lie groups. Their restricted root systems are their root systems of complex Lie algebra. So only $E_6$ remains. Category (b) consists of compact real forms of complex Lie algebras. But a compact real semisimple Lie group does not contain a proper parabolic subgroup. So its restricted root system is trivial. Complexifications of Lie algebras in category (c) are of $A$ to $D$ type. Root systems in category (d) is listed out in \cite[(6.108)]{Knapp2002}\footnote{There are two nonreduced root systems $BC_1,BC_2$ whose Weyl groups are $W(B_1)$ and $W(B_2)$.}. Only $E\mathrm{I}, E\mathrm{IV}$ remain.
\end{proof}
The following table\cite[(6.107)(6.108)]{Knapp2002} is all the real simple Lie algebras not excluded by Lemma~\ref{lemma:weyl_FPP}.
\begin{table}[htbp]
\centering
\label{tab:remaining-real-simple}
\begin{tabular}{lll}
\toprule
Real simple Lie algebra
& Cartan type
& Restricted root system \\
\midrule
$\mathfrak{sl}(n,\mathbb R)$, $n\geq 3$
& $A\mathrm{I}$& $A_{n-1}$ \\

$\mathfrak{sl}(n,\mathbb H)$, $n\geq 3$
& $A\mathrm{II}$ & $A_{n-1}$ \\

$\mathfrak{sl}(n,\mathbb C)_{\mathbb R}$, $n\geq 3$
& Complex type & $A_{n-1}$ \\

$\mathfrak{so}(n,n)$, $n\geq 3$ odd
& $D\mathrm{I}$ & $D_n$ \\

$\mathfrak{so}(2n,\mathbb C)_{\mathbb R}$,
$n\geq 3$ odd
& Complex type & $D_n$ \\

$E\mathrm{I}$
& $E\mathrm{I}$ & $E_6$ \\

$E\mathrm{IV}$
& $E\mathrm{IV}$ & $A_2$ \\

$E_6$
& Complex type & $E_6$ \\
\bottomrule
\end{tabular}
\end{table}

We now suppose that $P=P_I\subset G$ is defined by deleting $I$ from the Dynkin diagram of the restricted root system. Let $W_P$ be the Weyl group of this root system. By Bruhat decomposition, $G/P$ has cells indexed by the minimal-length coset representatives $W^P$ of $W/W_P$. Recall that rank of a compact Lie group is the dimension of its maximal torus.
\begin{lemma}\label{lemma:rank_FPP}
    Suppose that $K$ is a compact connected Lie group and $H$ is its closed subgroup. If $\rank H<\rank K$, then $K/H$ does not have the FPP.
\end{lemma}
\begin{proof}
    Let $T$ be the maximal torus of $K$ and $\overline{\langle t\rangle}= T$. Consider self-map on $f$ sending $kH$ to $tkH$. If $f(kH)=kH$, then $k^{-1}tk\in H$. Since $H$ is closed, $k^{-1}Tk\subset H$, which means $\rank H\ge \dim T=\rank K$.
\end{proof}
\begin{lemma}\label{lemma:parabolic_real_exc}
    $E\mathrm{I}/P$ does not have the FPP, if it is not isomorphic to $E\mathrm{I}/P_1$ or $E\mathrm{I}/P_6$. $E\mathrm{IV}/P$ does not have the FPP, if it is not isomorphic to $E\mathrm{IV}/P_1$ or $E\mathrm{IV}/P_2$.
\end{lemma}
\begin{proof}
    By Iwasawa decomposition, we have $E\mathrm{I}/P_I\cong K/Z_K(\mathfrak{a}_I)$.
    $\mathfrak{sp}(4)$ is the maximal compact subalgebra of $E\mathrm{I}$ \cite[Appendix~C.4]{Knapp2002}, so $\rank K=4$.  Denote that $L_I=Z_G(\mathfrak{a}_I)=Z_G(A_I)$ which is called the Levi factor of $G$. Let $S$ be the maximal torus of $Z_K(\mathfrak{a}_I)$ and its Lie algebra is $\mathfrak{s}$. By definition of $Z_K(\mathfrak{a}_I)$ and Cartan decomposition, $\mathfrak{s}\oplus\mathfrak{a}_I$ is an abelian Lie subalgebra of $E\mathrm{I}$. So
    \[
        \dim_\mathbb{R}\mathfrak{s}+\dim_\mathbb{R}\mathfrak{a}_I= \dim_\mathbb{C}(\mathfrak{s}\oplus\mathfrak{a}_I)_\mathbb{C}\le\rank \mathfrak{e}_6^\mathbb{C}=6.
    \]
    Since $E\mathrm{I}$ is a real split form and the restricted root system of a real split form is exactly the root system of its complexification \cite[Appendix C.4]{Knapp2002}, we obtain
    \[
        \rank Z_K(\mathfrak{a}_I)\le 6-|I|,
    \]
    with $\dim_\mathbb{R}(\mathfrak{a}_I)=|I|$. 
    This inequality only needs the split condition. It implies that if $|I|\ge 3$, $E\mathrm{I}/P_I$ does not have the FPP. The centralizer inverts the inclusion of deleted set, that is $Z_K(\mathfrak{a}_I)\subset Z_K(\mathfrak{a}_J)$ if $J\subset I$.
    If at least one of $2,3,4,5$ is contained in $I$, $\rank Z_K(\mathfrak{a}_I)\le 3$ since each $Z_K(\mathfrak{a}_i),i=2,3,4,5$ has rank $3$, see the table below \cite[Equation~(3.2)]{DobrevMarrani2020}.
    \begin{table}[ht]
\centering
\begin{tabular}{c|c|c}
\hline
$i$ & $\operatorname{Lie} Z_K(\mathfrak a_i)$
& $\operatorname{rank} Z_K(\mathfrak a_i)$ \\
\hline
$2$ & $\mathfrak{so}(6)$ & $3$ \\
$3$ or $5$ & $\mathfrak{so}(5)\oplus\mathfrak{so}(2)$ & $3$ \\
$4$ & $\mathfrak{so}(3)\oplus\mathfrak{so}(3)\oplus\mathfrak{so}(2)$ & $3$ \\
$1$ or $6$ & $\mathfrak{so}(5)\oplus\mathfrak{so}(5)$ & $4$ \\
\hline
\end{tabular}
\end{table}

Till now, remaining $I$ can only be $\{1\},\{6\},\{1,6\}$.
    \begin{center}        
% Preamble: \usepackage{tikz}

\begin{tikzpicture}[scale=1.15,
    every node/.style={font=\small}]

% Edges
\draw[thick] (0,0)--(4,0);
\draw[thick] (2,0)--(2,1);

% Dynkin nodes
\foreach \x in {0,4}{
    \filldraw[fill=gray!55,draw=black]
        (\x,0) circle (3.5pt);
}
\foreach \x in {1,2,3}{
    \filldraw[fill=red!65,draw=black]
        (\x,0) circle (3.5pt);
}
\filldraw[fill=red!65,draw=black]
    (2,1) circle (3.5pt);

% Bourbaki numbering (outside the nodes)
\node[below=6pt] at (0,0) {$1$};
\node[below=6pt] at (1,0) {$3$};
\node[below=6pt] at (2,0) {$4$};
\node[below=6pt] at (3,0) {$5$};
\node[below=6pt] at (4,0) {$6$};
\node[above=6pt] at (2,1) {$2$};
\node[below=6pt] at (2,-0.6) {$E_6$};

\end{tikzpicture}

    \end{center}
When $I=\{1,6\}$, let $w_0$ be the longest element in $W(E\mathrm{I},\mathfrak{a})$ which defined by 
\[
    w_0(\alpha_i)=-\alpha_{\sigma(i)},\quad \sigma=(16)(35).
\]
This element preserves $\alpha_{\{1,6\}}$ and switches $\alpha_1,\alpha_6$. Take the representative $\tau$ of $w_0$ and define self-map $f$ by sending $hZ_K(\mathfrak{a}_{\{1,6\}})$ to $h\tau Z_K(\mathfrak{a}_{\{1,6\}})$. This map does not have the FPP by the definition of $Z_K(\mathfrak{a}_{\{1,6\}})$.

For $E\mathrm{IV}$, the restricted root system is of type $A_2$. We only need to exclude $E\mathrm{IV}/P_{\{1,2\}}$ with minimal parabolic subgroup.
\begin{center}
    \begin{tikzpicture}[scale=1.15,
    every node/.style={font=\small}]

% Edges
\draw[thick] (0,0)--(1,0);

% Dynkin nodes
\foreach \x in {0,1}{
    \filldraw[fill=red!65,draw=black]
        (\x,0) circle (3.5pt);
}

% Numbering
\node[below=6pt] at (0,0) {$1$};
\node[below=6pt] at (1,0) {$2$};
\node[below=6pt] at (0.5,-0.6) {$A_2$};

\end{tikzpicture}
\end{center}
In this case, the Weyl group is $S_3$, any nontrivial element in $S_3$ defines an self-map by right multiplication which does not have the FPP.
\end{proof}

With the following two lemmas, these building blocks has the FPP by calculating the Lefschetz number of a self-map on their cohomology.
\begin{lemma}\label{lemma:EIV_exc}
    $E\mathrm{IV}/P_1\cong E\mathrm{IV}/P_2\cong\mathbb{O}P^2$.
\end{lemma}
\begin{proof}
    The first diffeomorphism comes from the symmetric of Dynkin diagram. The maximal compact subgroup of $E\mathrm{IV}$ is $F_4$ \cite[pp. 196]{RabeloSanMartin2019}. And the Levi subgroup is $\operatorname{Spin}(1,9)\times\mathbb{R}^+$ whose maximal compact subgroup is $\operatorname{Spin}(9)$ \cite[Table 1]{SugimotoShimokawa2022}. $F_4/\operatorname{Spin}(9)\cong\mathbb{O}P^2$ is proved in \cite[\S\,4.2]{RabeloSanMartin2019}
\end{proof}

\begin{lemma}\label{lemma:EI_exc}
    $E\mathrm{I}/P_1\cong E\mathrm{I}/P_6 \cong\Gr(2,\mathbb{H}^4)/\langle\tau\rangle$, where $\tau(V)=V^\perp$. And it has the same Betti numbers as $\mathbb{O}P^2$.
\end{lemma}
\begin{proof}
    The first diffeomorphism comes from the symmetric of Dynkin diagram.
    By Iwasawa decomposition, we have $E\mathrm{I}/P_1\cong K/Z_K(\mathfrak{a}_1)$.
    $\mathfrak{sp}(4)$ is the maximal compact subalgebra of $E\mathrm{I}$ \cite[Appendix~C.4]{Knapp2002}, so it is the Lie algebra of $K$. Consider $\widetilde{K}=\operatorname{Sp}(4)$ is a universal covering of $K$. Let $H=\pi^{-1}(Z_K(\mathfrak{a}_1))$ where $\pi$ is the covering map with $\ker\pi=\{\pm I_4\}$, then $E\mathrm{I}/P_1$ is diffeomorphic to $\operatorname{Sp}(4)/H$.

    Let $H^0$ be the component with identity of $H$. Since the semisimple part of the Levi subalgebra is $\mathfrak{so}(5,5)$ whose maximal compact subalgebra is $\mathfrak{so}(5)\oplus\mathfrak{so}(5)\cong \mathfrak{sp}(2)\oplus\mathfrak{sp}(2)$, Lie algebra of $H^0$ is $\mathfrak{sp}(2)\oplus\mathfrak{sp}(2)$. So $\rank H^0=\rank\widetilde{K}$ and $H^0$ and $\widetilde{K}$ have one common maximal torus. This induces an inclusion of restricted root systems $C_2\sqcup C_2\subset C_4$. Each inclusion determines a conjugation of $\operatorname{Sp}(2)\times\operatorname{Sp}(2)$, thus $\widetilde{K}/H^0\cong \operatorname{Sp}(4)/(\operatorname{Sp}(2)\times\operatorname{Sp}(2))\cong \Gr(2,\mathbb{H}^4)$.

    Now we suppose that $H^0=\operatorname{Sp}(2)\times\operatorname{Sp}(2)$ and set
    \[
        \tau=\begin{bmatrix}
            0 & I_2\\
            I_2 & 0
        \end{bmatrix}\in \operatorname{Sp}(4)\backslash H^0,
    \]
    which maps $V\in\Gr(2,\mathbb{H}^4)$ to $V^\perp$. And $N_{\widetilde{K}}(H^0)=H^0\sqcup H^0\tau $. Since $Z_G(\mathfrak{a}_1)=\operatorname{Spin}(5,5)\times \mathbb{R}^\times$ \cite[Table~1]{SugimotoShimokawa2022}, thus $Z_K(\mathfrak{a}_1)$ has two components. Since $\ker\pi\subset H^0$, $H$ has two components. $H^0$ is a normal subgroup of $H$, then $\tau\in H$. Therefore,
    \[
    \widetilde{K}/H\cong (\widetilde{K}/H^0)/(H/H^0)=\Gr(2,\mathbb{H}^4)/\langle\tau\rangle.
    \]

    Finally, the cohomology of the quaternionic Grassmannian is
    \[
        H^{4*}(\Gr(2,\mathbb{H}^4);\mathbb{Q})\cong\frac{\mathbb{Q}[q_1,q_2]}{(q_1^3-2q_1q_2,q_1^2q_2-q_2^2)}.
    \]
    The $\tau$ action switches between the tautological bundle $S$ and the quotient bundle $Q$. Since $S\oplus Q=\underline{\mathbb{H}}^4$, the symplectic Pontryagin classes satisfy $q(S)q(Q)=1$ and then
    \[
        \tau^*(q_1)=q_1(Q)=-q_1,\quad \tau^*(q_2)=q_2(Q)=q_1^2-q_2.
    \]
    This gives $H^*(E\mathrm{I}/P_1;\mathbb{Q})\cong\mathbb{Q}[q_1^2]/(q_1^6)$.
\end{proof}
\begin{remark}
     $H^*(E\mathrm{I}/P_1;\mathbb{Z})$ is not isomorphic to $\mathbb{Z}[q_1^2]/(q_1^6)$.
\end{remark}

Lemma~\ref{lemma:EIV_exc} and \ref{lemma:EI_exc} imply every self-map on $E\mathrm{I}/P_1$ or $E\mathrm{IV}/P_1$ induces an Adams map on the $\mathbb{Q}$-coefficient cohomology and gives the FPP.

The last building block is a complex simple Lie group. 
 Fix a maximal torus $T\subseteq B\subseteq P$, and let $W$ and $W_P$ denote the Weyl groups of $G$ and a Levi factor of $P$, respectively. By the Bruhat decomposition, $G/P$ admits a cellular decomposition into complex affine spaces indexed by the minimal-length coset representatives $ W^P$ of $W/W_P$, with each cell having complex dimension $\ell(w),w\in W^P$. Consequently, the integral cohomology of $G/P$ is torsion-free and concentrated in even degrees, with
\[
b_{2k}(G/P)=\#\{w\in W^P:\ell(w)=k\}.
\]
In particular, its half-degree Poincaré polynomial is given by
\[
\tilde{P}_{G/P}(t)
=\sum_{w\in W^P}t^{\ell(w)}
=\frac{\sum_{w\in W}t^{\ell(w)}}{\sum_{w\in W_P}t^{\ell(w)}},
\]
and its Euler characteristic is $\chi(G/P)=|W|/|W_P|$. Moreover, the rational cohomology ring admits Borel's presentation
\[
H^*(G/P;\mathbb Q)
\cong
\frac{\operatorname{Sym}(\mathfrak t_{\mathbb Q}^{*})^{W_P}}
{\bigl(\operatorname{Sym}(\mathfrak t_{\mathbb Q}^{*})_{>0}^{W}\bigr)},
\]
where $\mathfrak t_{\mathbb Q}^{*}$ consists of weights whose elements are assigned cohomological degree $2$. See \cite{BGG1973} for more details. With these setups, we have the following statement on complex simple Lie groups.
\begin{proposition}\label{prop:E_6_exc}
    Suppose that $G$ is a connected complex simple Lie group of exceptional type and $P$ is a proper parabolic subgroup. Then $G/P$ has the FPP if and only if it is isomorphic to $E_6/P_1$.
\end{proposition}
\begin{proof}
    Let $X_0$ be the split real form of $G/P$ induced by the split real form of $G$ with a complex conjugation $\sigma$ on it. Moreover,
    \[
        \sigma^*\big|_{H^{2j}(G/P;\mathbb{Q})}=(-1)^j\operatorname{Id}.
    \]
    and then the Leftschetz number
    \[
        L(\sigma)=\chi(X_0) = \tilde{P}_{G/P}(-1) = \left.\frac{\sum_{w\in W}(t)^{\ell(w)}}{\sum_{w\in W_P}(t)^{\ell(w)}}\right|_{t=-1},
    \]
    which is nonzero if and only if $G/P$ is $E_6/P_1$, $E_6/P_6$ and $E_6/P_{\{1,6\}}$, see Appendix~\ref{sec:lef_num}. Since $G/P$ is simply connected, then the $L(\sigma)=0$ implies that Neilsen number is zero by the main theorem in \cite{Jiang1983}.
    Then by Wecken's Theorem \cite[Satz~3]{Wecken1942} \cite[Theorem~2.1]{Jezierski2003}, $\sigma$ is homotopic to a map with no fixed point. Then the only possible $G/P$ who has the FPP in these three cases.

    Since the longest Weyl group element $w_0$ of $E_6$ preserves the root subsystem of type $D_4$ corresponding to $\{2,3,4,5\}$, a representative $a\in N_K(T)$ normalizes $U=K\cap P_{\{1,6\}}$, where $K$ is the maximal compact subgroup. Moreover, $w_0\notin W_{P_{\{1,6\}}}$, so $a\notin U$. Hence right multiplication by $a$ induces a fixed-point-free diffeomorphism of $K/U\cong E_6/P_{\{1,6\}}$. Consequently, $E_6/P_{\{1,6\}}$ does not have the fixed-point property. The following is Dynkin diagram of $E_6/P_{\{1,6\}}$.
    \begin{center}
        
% Preamble: \usepackage{tikz}

\begin{tikzpicture}[
    scale=1.15,
    every node/.style={font=\small}
]

% Solid edges
\draw[thick] (1,0)--(3,0);
\draw[thick] (2,0)--(2,1);

% Dashed edges adjacent to nodes 1 and 6
\draw[thick, dashed] (0,0)--(1,0);
\draw[thick, dashed] (3,0)--(4,0);

% Nodes 1 and 6: dashed circles
\foreach \x in {0,4}{
    \draw[thick, dashed, fill=gray!20]
        (\x,0) circle (4pt);
}

% Nodes 3, 4, 5: solid circles
\foreach \x in {1,2,3}{
    \draw[thick, fill=red!65]
        (\x,0) circle (4pt);
}

% Node 2: solid circle
\draw[thick, fill=red!65]
    (2,1) circle (4pt);

% Bourbaki numbering
\node[below=7pt] at (0,0) {$1$};
\node[below=7pt] at (1,0) {$3$};
\node[below=7pt] at (2,0) {$4$};
\node[below=7pt] at (3,0) {$5$};
\node[below=7pt] at (4,0) {$6$};
\node[above=7pt] at (2,1) {$2$};
\node[below=6pt] at (2,-0.6) {$E_6$};

\end{tikzpicture}

    \end{center}

By \cite[Theorem~3]{LiDuan2009}, every graded endomorphism of $H^*(E_6/P_1;\mathbb{Z})$ is an Adams map\footnote{Their $E_6$ is the compact exceptional Lie group and $E_6/(D_5\cdot S^1)$ is our $E_6/P_1$}. Then for every self-map $f$ on $E_6/P_1$ with $f^*|_{H^{2j}}=m^j\operatorname{Id}$, we have
\[
    L(f) = \frac{P_{W(E_6)}(m)}{P_{W(D_5)}(m)}=[9]_m[3]_{m^4}>0,\quad m\in\mathbb{Z},
\]
where $[k]_m=1+m+\cdots+m^{k-1}$.
So $E_6/P_1\cong E_6/P_6$ has the FPP.
\end{proof}

\subsection*{Proof of Theorem~\ref{mainthm:real_exc_FPP}} 
Suppsoe that $G/P$ has the FPP, then each simple component has the FPP. By Corollary~\ref{cor:exc_remaining}, Lemma~\ref{lemma:parabolic_real_exc}, \ref{lemma:EIV_exc}, \ref{lemma:EI_exc} and Proposition~\ref{prop:E_6_exc}, each simple component must be in three types.

Conversely, we write $X=\mathbb{O}P^2$, $Y=\Gr(2,\mathbb{H}^4)/\langle\tau\rangle$ and $Z=E_6/P_1$. We first prove that $X^a\times Y^b$ has the FPP for some $a,b$. These two flag manifolds have the same rational coefficient cohomology. By K\"{u}nneth formula,
\[
    H^{8*}(X^a\times Y^b;\mathbb{Q})=\frac{\mathbb{Q}[u_1,\ldots,u_r]}{(u_1^3,\ldots,u_r^3)},\quad r=a+b.
\]
For any self-map $f$, suppose that
\[
    f^*(u_i)=\sum_j a_{ij}u_j.
\]
The relation gives
\[
    0=f^*(u_i^3)=\left(\sum_j a_{ij}u_j\right)^3,
\]
which means $a_{ij}^2a_{ik}=0, k\ne j$ and matrix $[a_{ij}]$ has at most one nonzero entry for each row. Thus, $f^*(u_i)=a_iu_{\sigma(i)}$ for some self-map $\sigma$ on $\{1,\ldots,r\}$. It defines a directed graph with vertex $i$ and edge $i\to \sigma(i)$. If $C$ is one cycle of the graph, define $a_C$ to be the product of those $a_i$ with $i\in C$. We calculate the trace with respect to the monomial basis $u_1^{e_1}\cdots u_r^{e_r},0\le e_i\le 2$. A monomial basis contributes to the diagonal if and only if it is supported by cycles and $e_i$ are the same in each cycle. Denote this diagonal element by a product of $a_C^{e_C}$.
\[
    L(f) =\sum_{e_{C_i}}\prod_j a_{C_j}^{e_{C_j}}=\prod_C (1+a_C+a_C^2)>0,
\]
where $C$ runs through all circles in the graph.

Then we prove the theorem for $Z^c$ for some $c$. Consider $\mathbb{F}_2$-coefficient cohomology \cite[Equation (2.5)]{LiDuan2009},
\[
    A:=H^*(Z;\mathbb{F}_2)\cong \frac{\mathbb{F}_2[h,s]}{(hs^2,s^3+h^{12})},\quad \deg h = 2,\deg s =8.
\]
Since odd degrees of $A$ vanish, we have a complex $\operatorname{Sq}^2:H^i\to H^{i+2}$ by Steenrod square. Its cohomology is 
\[
 H^*(A,\operatorname{Sq}^2)\cong \mathbb{F}_2[v]/(v^3), \quad v=[s^2].
\]
and we can calculate the Lefschetz number on this complex since the natruality of Steenrod square $f^*\operatorname{Sq}^2=\operatorname{Sq}^2f^*$. By K\"{u}nneth formula, $H^*(Z^c;\mathbb{F}_2)\cong A^{\otimes c}$ and
\[
    H^*(A^{\otimes c},d)\cong\frac{\mathbb{F}_2[v_1,\ldots,v_c]}{(v_1^3,\ldots,v_c^3)},\quad d=\sum_i 1^{\otimes(i-1)}\operatorname{Sq}^2\otimes 1^{\otimes (c-i)}.
\]
By the argument in $X^a\times Y^b$ case, the Lefschetz number $L(f)\equiv 1 \mod 2$ 
on the complex $(A^{\otimes c},d)$, then $Z^c$ has the FPP.

Finally, we prove the mixed case $X^a\times Y^b\times Z^c$. Denote that $B=H^*(X^a\times Y^b;\mathbb{Q})$ and $C=H^*(Z^c;\mathbb{Q})$.

\noindent \textbf{Claim.} There is no nontrivial $x\in H^8(Z^c;\mathbb{Q})$ satisfies $x^3= 0$. We first verified it for $\mathbb{F}_2$-coefficient cohomology.
When $c=1$, 
\[
    H^8(Z;\mathbb{F}_2)\cong \Span_{\mathbb{F}_2}(h_i^4,s_i:i=1,\ldots,c).
\]
it is easy to verify by simple calculations. When $c>1$, write
\[
    x = \sum_i a_i h_i^4 +b_i s_i +Q(h_1,\ldots,h_c),\quad a_i,b_i\in\mathbb{F}_2,
\]
where $Q$ is a homogeneous polynomial of degree $4$ and does not contain any $h_i^4$. By argument for $c=1$, if $x^3=0$, we have $a_i=b_i=0$ and $Q^3=0$. By K\"{u}nneth formual monomials of $h_i$ in degree $12$ are linear independent. So $Q=0$. Since $H^*(Z^c;\mathbb{Z})$ is torsion-free, the claim holds for rational coefficient cohomology.

Then we prove that any self map $f$ on $X^a\times Y^b\times Z^c$ preserves $B$. Suppose that
\[
    f^*(u_i) = \sum_j a_{ij}u_j + x_i,\quad 1\le i\le a+b, \quad x_i\in H^8(Z^c;\mathbb{Q}).
\]
Since $f^*(u_i^3)=0$, each $x_i=0$ by the claim. So $f^*(B)\subset B$. Define the following two self-maps
\[
    f_1: m\mapsto p_{X\times Y}\circ f(m,z_0),\quad f_2:z\mapsto p_{Z}\circ f(m_0,z),\quad m,m_0\in X^a\times Y^b,\quad z,z_0\in Z^c.
\]
where $p_{X\times Y}$ and $p_Z$ are projections to $X^a\times Y^b$ and $Z^c$ respectively. 
We define
\[
    R=f^*-f^*_1\otimes f^*_2: B\otimes C\to B\otimes C.
\]
Since $f^*(B)\subset B$, for every $b\in B,c\in C$, $f^*(b\otimes 1)=f_1^*(b)$ and
\[
    R(b\otimes c) = f^*(b\otimes c)-f^*_1(b)\otimes f^*_2(c)= (f_1^*(b)\otimes1)(f^*(1\otimes c)-1\otimes f_2^*(c)).
\]
By the definition of $f_2$, 
\[
    f^*(1\otimes c)-1\otimes f_2^*(c)\in B_{>0}\otimes C.
\]
Since $f_1^*$ keeps the degree, $R$ increases the degree on component $B$, which means $R$ is a nilpotent map. So trace of $\operatorname{tr}(f^*)=\operatorname{tr}(f_1^*)\operatorname{tr}(f_2^*)$. Since the cohomology of the space has no odd degree, $L(f)=\operatorname{tr}(f^*)$.
With the argument in the former two cases, $L(f)\ne 0$. \hfill $\Box$

\subsection*{Proof of Theorem~\ref{mainthm:BCD_FPP}}
By Corollary~\ref{cor:BCD_FPP}, we only need to prove it for a real simple Lie group $G$ with $D_{n}, n=2m+1\ge 5$ type restricted root system. There are only two cases whose Lie algebras are $\mathfrak{so}(2m+1,2m+1)$ and $\mathfrak{so}(4m+2,\mathbb{C})_\mathbb{R}$.

For $G=\operatorname{SO}_0(n,n)$, consider a quadratic form $Q(x,y)=\|x\|^2-\|y\|^2$ on $\mathbb{R}^n\oplus\mathbb{R}^n$. Every flag manifold can be defined by a partial flag of isotropic subspaces.

If this defining flag does not contain a maximal isotropic subspace, let
\[
    T=\begin{bmatrix}
        -I_n & 0\\ 0 & I_n
    \end{bmatrix}\in\operatorname{O}(n,n),
\]
which preserves $Q$ and then the $G$-orbit of defining isotopic flag. So it induces a self-map on $G/P$ such that on each isotropic flag 
\[
    f_T:V_1\subset \cdots \subset V_k\mapsto TV_1\subset \cdots \subset TV_k.
\]
If there is a $W$ such that $TW=W$. Every $(x,y)\in W$ satisfies $(-x,y)=T(x,y)\in W$. Then $(x,0),(0,y)\in W$ and $Q(x,0)=Q(0,y)=0$. Thus, $(x,y)=0$. This means that $f$ has no fixed point.

If the defining flag contains a maximal isotropic subspace, then every isotropic flags in $G/P$ contains a maximal isotropic subspace. Let $A=\operatorname{diag}(-1,-1,1,\ldots,1)\in \operatorname{SO}(n)$ and
\[
    S=\begin{bmatrix}
        I_n & 0\\ 0 & A
    \end{bmatrix}\in\operatorname{SO}_0(n,n),
\]
which induces a self-map $f_S$ on $G/P$. For a maximal isotropic subspace $W$, it is isomorphic to the graph $\Gamma_B=\{(x,Bx):x\in\mathbb{R}^n\}$ for some $B\in\operatorname{O}(n)$. $S(W)$ is isomorphic to $\Gamma_{AB}$. If $S(W)=W$, then $AB=B$ by the uniqueness of graph. Then $A=I_n$ since $B$ is invertible, which is a contradiction. So $f_S$ has no fixed point.

For $G=\operatorname{SO}(2n,\mathbb{C})$, consider $V=\mathbb{C}^{2n}$ and $Q=z_1^2+\cdots+z_{2n}^2$. 

If the defining flag does not contain a maximal isotropic subspace, then the complex conjugation $\sigma$ keeps the isotropic condition with respect to $Q$ which induces a self-map $f_\sigma$ on $G/P$ by sending flags to their conjugations. If $\overline{W}=W$ for some isotropic subspace, then every $z=x+\sqrt{-1}y\in W$ satisfies $x,y\in W$. Then $Q(x)=Q(y)=0$ gives $x=y=0$. So $W=0$ and $f_\sigma$ has no fixed point.

If the defining flag contains a maximal isotropic space, we define
\[
    \tau(z_1,\ldots,z_{2n})=(\overline{z_1},\ldots,\overline{z_{2n-1}},-\overline{z_{2n}}).
\]
It keeps the isotropic condition with respect to $Q$ and induces a self-map $f_\tau$. If $W$ is a maximal isotropic subspace such that $\tau(W)=W$, there must be a nonzero $z=(z_1,\ldots,z_{2n-1},0)\in W$ since its complex dimension is $n$. Suppose that $z=x+\sqrt{-1}y$, then $z,\overline{z}\in W$ gives $x,y\in W$ and $Q(x)=Q(y)=0$. Thus, $z=0$ which is a contradiction. So $f_\tau$ has no fixed point.\hfill$\Box$

\appendix

\section{Lefschetz numbers on homogeneous spaces}\label{sec:lef_num}

Suppose that $G$ is a connected complex simple  Lie group of type $E_6$ and $P_I$ is the parabolic subgroup defined by deleting $I$ from its Dynkin diagram. Write $G/P_I\cong X_0\times_{\Spec\mathbb{R}}\Spec\mathbb{C}$ and $\sigma$ is the conjugation. We calculate the $L(\sigma)$ when $I\ne\{1\},\{6\}$ or $\{1,6\}$.

The half-degree Poincar\'{e} polynomial is
\[
    \tilde{P}_{G/P_I}(t)= \frac{\tilde{P}_{W(E_6)}(t)}{\tilde{P}_{W_I}(t)}=\frac{[2]_t[5]_t[6]_t[8]_t[9]_t[12]_t}{\tilde{P}_{W_I}(t)},
\]
where $t+1$ has multiple $4$ in the enumerator. $L(\sigma)=\tilde{P}_{G/P_I}(-1)$. Use the following table, we obtain that only on $E_6/P_1$, $E_6/P_6$ or $E_6/P_{\{1,6\}}$, $\sigma$ has a nonzero Lefschetz number.
\begingroup
\footnotesize
\setlength{\tabcolsep}{5pt}
\renewcommand{\arraystretch}{1.2}

\begin{longtable}{
@{}p{0.53\textwidth}
p{0.29\textwidth}
c@{}
}

\toprule
Deleted nodes $I$
& $\tilde{P}_{W_I}(t)$
& $\operatorname{ord}_{t=-1}$\\
\midrule
\endfirsthead

\toprule
Deleted nodes $I$
& $P_{W_I}(t)$
& $\operatorname{ord}_{t=-1}P_{W_I}(t)$\\
\midrule
\endhead

$1,6$
& $[2][4][5][6][8]$ & $4$\\

$2$
& $[2][3][4][5][6]$ & $3$\\

$3,5$
& $[2]^2[3][4][5]$ & $3$\\

$4$
& $[2]^3[3]^2$ & $3$\\

$16$
& $[2][4]^2[6]$ & $4$\\

$12,13,26,56$
& $[2][3][4][5]$ & $2$\\

$15,23,25,36$
& $[2]^2[3][4]$ & $3$\\

$24$
& $[2]^2[3]^2$ & $2$\\

$14,34,35,45,46$
& $[2]^3[3]$ & $3$\\

$123,126,136,156,256$
& $[2][3][4]$ & $2$\\

$124,125,134,135,234,$
\newline
$236,245,246,356,456$
& $[2]^2[3]$ & $2$\\

$145,146,235,345,346$
& $[2]^3$ & $3$\\

$1234,1236,1256,1356,2456$
& $[2][3]$ & $1$\\

$1235,1245,1246,1345,1346,$
\newline
$1456,2345,2346,2356,3456$
& $[2]^2$ & $2$\\

$12345,12346,12356,$
\newline
$12456,13456,23456$
& $[2]$ & $1$\\

$123456$
& $1$ & $0$\\

\bottomrule
\end{longtable}
\endgroup

\section{Fixed Point Property of Flag mainfolds}\label{sec:append_flag_FPP}

\begin{longtable}{@{}llll@{}}
\toprule
Type & Real Lie algebra & Restricted roots & Flag manifolds with FPP \\
\midrule
\endfirsthead

\toprule
Type & Real Lie algebra & Restricted roots & Flag manifolds with FPP \\
\midrule
\endhead

\bottomrule
\endlastfoot

\multicolumn{4}{@{}l}{\textbf{I. Compact simple Lie groups}}\\
\addlinespace[2pt]
& All compact simple Lie algebras
& $\varnothing$
& None \\

\midrule

\multicolumn{4}{@{}l}{
    \textbf{II. Complex simple Lie groups regarded as real Lie groups}
}\\
\addlinespace[2pt]

$A_{n-1}$ & $\mathfrak{sl}(n,\mathbb{C})$
    & $A_{n-1}$
    & \begin{tabular}[t]{@{}l@{}}
        $\Gr(k,\mathbb{C}^N)$: classified\\
        General: conjectural
      \end{tabular} \\
$B_n$ & $\mathfrak{so}(2n+1,\mathbb{C})$
    & $B_n$ & None \\
$C_n$ & $\mathfrak{sp}(2n,\mathbb{C})$
    & $C_n$ & None \\
$D_n$ & $\mathfrak{so}(2n,\mathbb{C})$
    & $D_n$ & None \\
$E_6$ & $\mathfrak{e}_6(\mathbb{C})$
    & $E_6$ & $E_6/P_1 \cong E_6/P_6$ \\
$E_7$ & $\mathfrak{e}_7(\mathbb{C})$
    & $E_7$ & None \\
$E_8$ & $\mathfrak{e}_8(\mathbb{C})$
    & $E_8$ & None \\
$F_4$ & $\mathfrak{f}_4(\mathbb{C})$
    & $F_4$ & None \\
$G_2$ & $\mathfrak{g}_2(\mathbb{C})$
    & $G_2$ & None \\

\midrule

\multicolumn{4}{@{}l}{
    \textbf{III. Noncompact classical real forms}
}\\
\addlinespace[2pt]

AI & $\mathfrak{sl}(n,\mathbb{R})$
    & $A_{n-1}$
    & \begin{tabular}[t]{@{}l@{}}
        $\Gr(k,\mathbb{R}^N)$: classified\\
        General: conjectural
      \end{tabular} \\

AII & $\mathfrak{su}^{*}(2n)$
    & $A_{n-1}$
    & \begin{tabular}[t]{@{}l@{}}
        $\Gr(k,\mathbb{H}^N)$: classified\\
        General: conjectural
      \end{tabular} \\
AIII & $\mathfrak{su}(p,q)$, $p<q$
    & $BC_p$ & None \\
AIII & $\mathfrak{su}(p,p)$
    & $C_p$ & None \\
BI, DI & $\mathfrak{so}(p,q)$, $p<q$
    & $B_p$ & None \\
DI & $\mathfrak{so}(p,p)$
    & $D_p$ & None \\
CI & $\mathfrak{sp}(2n,\mathbb{R})$
    & $C_n$ & None \\
CII & $\mathfrak{sp}(p,q)$, $p<q$
    & $BC_p$ & None \\
CII & $\mathfrak{sp}(p,p)$
    & $C_p$ & None \\
DIII & $\mathfrak{so}^{*}(4n)$
    & $C_n$ & None \\
DIII & $\mathfrak{so}^{*}(4n+2)$
    & $BC_n$ & None \\

\midrule

\multicolumn{4}{@{}l}{
    \textbf{IV. Noncompact exceptional real forms}
}\\
\addlinespace[2pt]

EI & $\mathfrak{e}_{6(6)}$
    & $E_6$ & $\mathrm{EI}/P_1 \cong \mathrm{EI}/P_6$ \\
EII & $\mathfrak{e}_{6(2)}$
    & $F_4$ & None \\
EIII & $\mathfrak{e}_{6(-14)}$
    & $BC_2$ & None \\
EIV & $\mathfrak{e}_{6(-26)}$
    & $A_2$ & $\mathrm{EIV}/P_1 \cong \mathrm{EIV}/P_2$ \\
EV & $\mathfrak{e}_{7(7)}$
    & $E_7$ & None \\
EVI & $\mathfrak{e}_{7(-5)}$
    & $F_4$ & None \\
EVII & $\mathfrak{e}_{7(-25)}$
    & $C_3$ & None \\
EVIII & $\mathfrak{e}_{8(8)}$
    & $E_8$ & None \\
EIX & $\mathfrak{e}_{8(-24)}$
    & $F_4$ & None \\
FI & $\mathfrak{f}_{4(4)}$
    & $F_4$ & None \\
FII & $\mathfrak{f}_{4(-20)}$
    & $BC_1$ & None \\
G & $\mathfrak{g}_{2(2)}$
    & $G_2$ & None \\

\end{longtable}

\printbibliography

\end{document}